\documentclass[journal]{aiaa-pretty}
\usepackage[utf8]{inputenc}
\usepackage{textcomp}

\usepackage{graphicx}
\usepackage{amsmath}
\usepackage[version=4]{mhchem}
\usepackage{siunitx}
\usepackage{longtable,tabularx}
\usepackage{subfigure}
\usepackage{cite}

\usepackage{booktabs}
 \usepackage{varioref}
 \usepackage{wrapfig}
 \usepackage{threeparttable}
 \usepackage{dcolumn}
  \newcolumntype{d}{D{.}{.}{-1}}

\usepackage{letltxmacro}
\LetLtxMacro{\originaleqref}{\eqref}
\renewcommand{\eqref}{Eq.~\originaleqref}

\newcommand{\bC}{{\boldsymbol C}}
\newcommand{\bF}{{\boldsymbol F}}

\newcommand{\bi}{{\boldsymbol i}}
\newcommand{\bj}{{\boldsymbol j}}
\newcommand{\bk}{{\boldsymbol k}}

\newcommand{\br}{{\boldsymbol r}}

\newcommand{\bu}{{\boldsymbol u}}

\newcommand{\bA}{{\boldsymbol A}}
\newcommand{\bB}{{\boldsymbol B}}

\newcommand{\bT}{{\boldsymbol T}}

\newcommand{\bomega}{\mbox{\boldmath$\omega$}}

\newcommand{\Real}{\mathbb R}
\newcommand{\set}[1]{\left\{#1\right\}}
\newcommand{\real}[1]{{\mathbb R}^{#1}}
\title{Nechvile-Transformed Spacecraft Dynamics and Propellant Computation in the 3-Body Problem\footnote{A version of this paper was published in the Jan 2026 issue of the \textit{Journal of Guidance Control and Dynamics} under the heading, ``Nuances in Propellant Computation in the Elliptic Restricted Three-Body Problem;''  see https://doi.org/10.2514/1.G009053  }}

\author{Michael J. Dixon\footnote{CDR, USN, and Graduate Student, Department of Mechanical and Aerospace Engineering, Monterey, CA. michael.dixon@nps.edu} and Isaac M. Ross\footnote{Distinguished Professor, Department of Mechanical and Aerospace Engineering, Monterey, CA.}
\\
\textit{Naval Postgraduate School, Monterey, CA 93943}
}

\abstract{
The uncontrolled equations of motion in the Nechvile frame for the restricted three-body problem have been well-known since at least the 1960s. It would seem that adding an external force to these equations is quite trivial: simply add an external force per mass term to the acceleration equations. Here we show that the last statement is not true. In fact, we show that the additive generic external force must be multiplied by the inverse of $(1 + e \cos\theta)^3$ where $e$ is the relative eccentricity of the primaries and $\theta$ is the true anomaly of the rotating frame located at the barycenter. Furthermore, when this result is combined with the mass flow rate equation, it generates several surprising results due to the mismatch between the resulting quadratic term and the cubic term in the equations of motion. This leads to a corresponding modification of the rocket equation itself.  A Birkhoff-theoretic solution to an illustrative cislunar space mission problem shows propellent savings of 80\% with the use of the correct cost functional.   The popular quadratic cost utilizes more than $2X$ the minimum propellant consumption.
}

\begin{document}

\maketitle

\section{Introduction}\label{sec:intro}
In principle, computing propellant consumption for any given orbital maneuver is quite easy: simply integrate the rocket mass flow rate equation over the thrust profile\cite{biblarz,hale}.  If the thrust profile is given in terms of a timed sequence of impulsive vectors, then the integral of the mass flow rate equation reduces to a sum of the famous Tsiolkovsky formula\cite{hale}.  These foundational principles are well-understood. Hence, at first glance, there does not seem to be any reason to revisit them for spacecraft motion in the elliptic restricted three-body problem (ER3BP). Except that this paper proves otherwise for the most widely-used\cite{szebehely_theory_1967,broucke_stability_1969,celletti_dynamics_2024, du_low-thrust_2023,quartullo-2023, peng-2015} equations of motion for the ER3BP.

The equations of motion for the ER3BP are given most elegantly through the use of Nechvile’s transformation\cite{szebehely_theory_1967,broucke-nechville-1966, peng-2015,duboshin-1971}. In addition to changing the independent variable from time to the true anomaly, Nechvile’s transformation renders the rotating frame to be ``pulsating'' and maintains the relative stationarity of the Lagrange points\cite{szebehely_theory_1967,peng-2015}.  The price for this elegance is that the thrust force in the equations of motion must be multiplied by the inverse of a cubic term given by\cite{du_low-thrust_2023,quartullo-2023},
\begin{equation}\label{eq:cubic}
(1+e\cos\theta)^3
\end{equation}
where $e$ is the relative eccentricity of the two primaries and $\theta$ is their true anomaly. In this paper, we show that an important fallout of Nechvile's transformation is that the rocket mass flow rate equation must also be modified but by a different multiplicative term given by the inverse of the quadratic,
\begin{equation}\label{eq:quadratic}
(1+e\cos\theta)^2
\end{equation}
Because the term given by \eqref{eq:quadratic} is not a constant and different from \eqref{eq:cubic}, we show that the computation of propellant consumption in the ER3BP must be modified in a new manner.  Obviously, the terms in both Eqs.~(\ref{eq:cubic}) and (\ref{eq:quadratic}) reduce to unity if $e=0$; i.e., the circular restricted three-body problem (CR3BP). In other words, the CR3BP  ``hides'' the effects of  Eqs.~(\ref{eq:cubic}) and (\ref{eq:quadratic}).

In the early stages of an analysis of a space mission, it is desirable to eliminate the mass flow rate equation while simultaneously minimizing propellant consumption.  Such mission designs are agnostic to the specifics of a rocket propulsion system. This is part of the origin of the impulsive Delta-V maneuvers\cite{robbins-deltaV}. Because practical rockets are not impulsive, the concept of minimizing gravity loss and drag loss (where applicable) took hold\cite{biblarz,hale,robbins-deltaV}. It was shown in [\citenum{ross_how_2004,ross_space_2006}] that this two-step optimization process can be reduced to a single integrated step through the use of an $L^1$ cost functional. Although the $L^1$ functional is not smooth, it can be easily mapped to a smooth differentiable function through the use of elementary mathematical transformations\cite{ross-book}.  These ideas obviate the need for using the popular but erroneous quadratic cost fucntionals\cite{ross_how_2004,ross_space_2006,ross-book} as proxies for propellant consumption.  In any event, because of the
``discrepancy'' between Eqs.~(\ref{eq:cubic}) and (\ref{eq:quadratic}), we show that the $L^1$ functional proposed in [\citenum{ross_how_2004, ross_space_2006}] gets modified when used in conjunction with the Nechvile-transformed equations of motion. In addition, we show that even the classic \mbox{Delta-V} computation must also be modified by a linear cosine term given by $(1 + e \cos\theta)$.

If the cubic term given by \eqref{eq:cubic} is used to simply scale the thrust force using Nechvile's units, then the control space must be treated as time-varying in order to account for the full capability of a practical rocket engine\cite{biblarz}. If this key point is ignored and the control space is simply boxed for mathematical convenience\cite{chyba,caillau}, then it is shown that the utilization of a practical engine is significantly reduced even for small eccentricities of the primaries. For instance, in the case of a spacecraft in cislunar space ($e \simeq 0.05$), the artificial constraining of the control space limits the usage of a practical rocket to less than $75\%$ of its maximum capability.

If the preceding nuances in analysis and computation are ignored for a practical space mission, the best case scenario is an increase in propellant consumption. The worst case possibility is a mission failure. These findings and more are numerically illustrated for a sample cislunar mission in the ER3BP.

\section{Equations of Motion for the Third Body Subject to an External Force}

In the restricted three-body problem\cite{szebehely_theory_1967}, the mass of the third body is negligible compared to those of the two primaries.  We assume the third body is a spacecraft of mass $m_{sc}$ and the two primaries are of masses $m_1$ and $m_2$. Hence, we have
$$ m_1 +m_2 + m_{sc} \approx m_1 + m_2$$
We use the standard mass parameter, $\mu$, defined by\cite{szebehely_theory_1967,broucke_stability_1969},
\begin{equation}
    \mu:=\frac{m_2}{m_1+m_2}
    \label{eq: massParameter_mu}
\end{equation}
By convention $m_2 \le m_1$; hence, $\mu \in (0, 0.5]$.

In a remarkable paper~\cite{elipe-wrongHyp-2024}, Elipe proves that a large number of results in the restricted three-body problem ``published even in well-reputed journals ... are wrong because their hypotheses are against physical laws.''  Consequently, it is necessary to exercise great caution in deriving/modifying equations of motion for the restricted three-body problem under new/added hypotheses. The uncontrolled equations of motion for the ER3BP are well-known\cite{szebehely_theory_1967}; however, incorporating an external force to these uncontrolled equations is more nuanced than simply adding control terms.  In following Elipe's caution\cite{elipe-wrongHyp-2024} in using ``off-the-shelf'' equations of motion from the literature, we start with first principles using Szebehely's derivation\cite{szebehely_theory_1967} to produce the same for a ``spacecraft'' in the ER3BP. The extra hypothesis we make herein is that the third body (i.e., spacecraft) is subject to an arbitrary external force. This external force does not act on the other two bodies.

Szebehely's original work is followed where appropriate, but the methodologies in this paper utilize standard techniques from dynamics.  For example, instead of using Birkhoff's complex variable transformation\cite{birkhoff-1915} as done by Szebehely, we use standard coordinate transformations from the inertial frame to the rotating frame as in \cite{gurfil_stationkeeping_2012}.  An alternative derivation using Lagrangian mechanics is provided in [\citenum{quartullo-2023}].

\subsection{Variables in the ER3BP}\label{sec:basic}

\begin{figure}[hbt!]
\centering
\subfigure[Schematic for a spacecraft in the ER3BP]{%
\centering
    \includegraphics[width=\columnwidth, trim={1.75in 0.5in 5.5in 1in},clip]{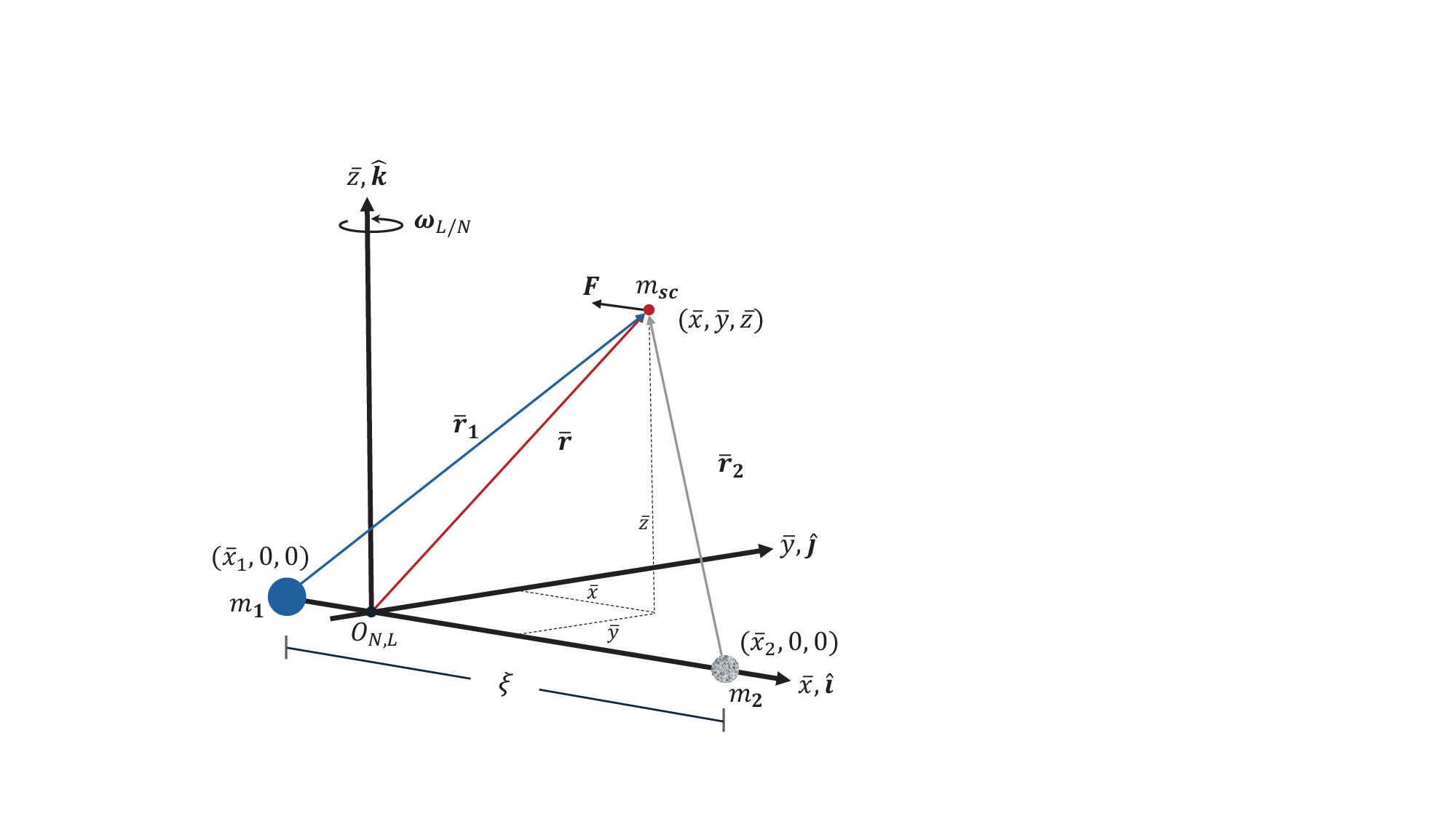}
}%
\hfill
\subfigure[Planar view of the rotating barycentric frame]{%
\centering
    \includegraphics[width=0.8\columnwidth, trim={1in 0in 0in 0in},clip]{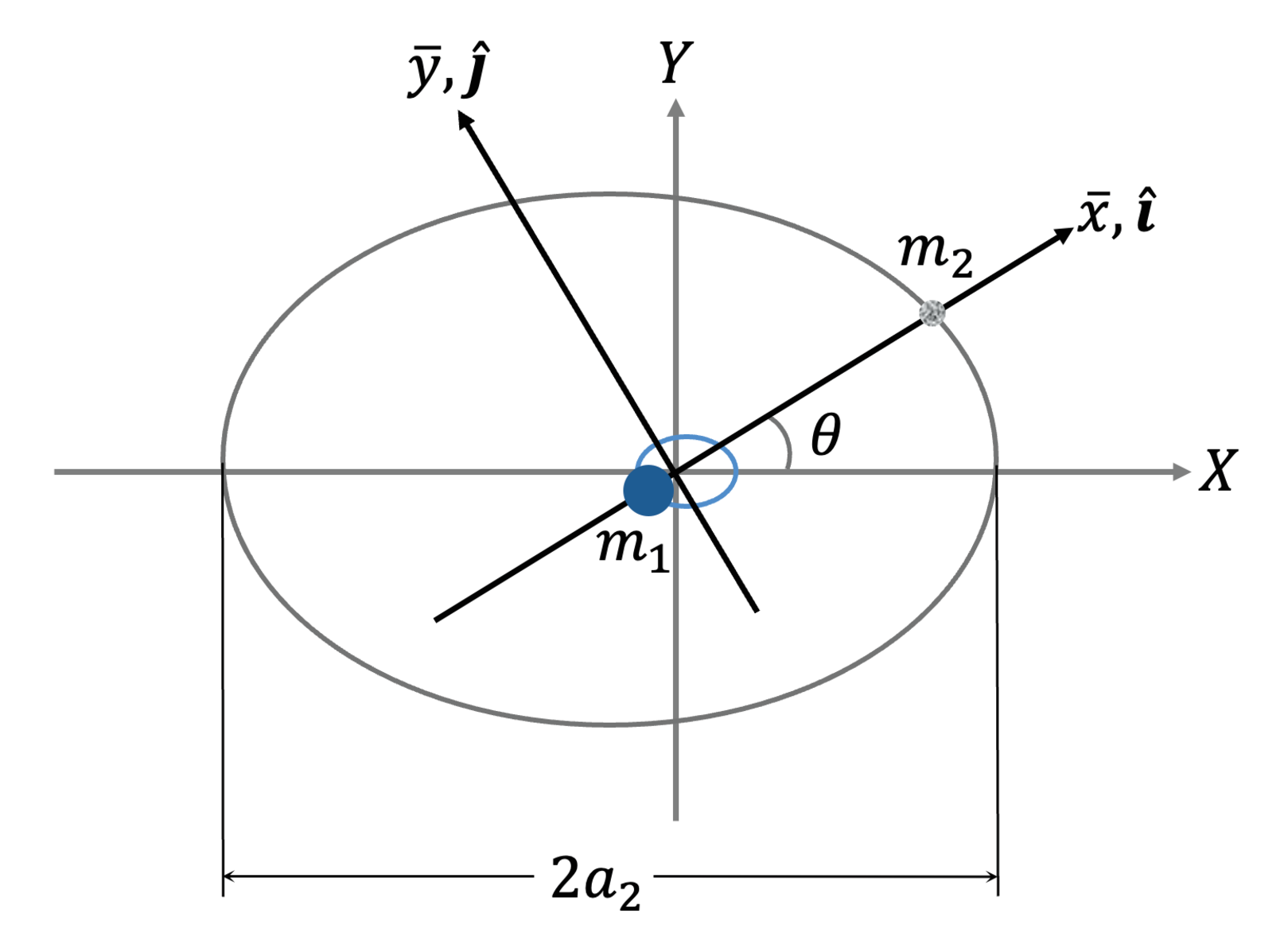}
    \label{fig: CoordFrame_BaryToInertial_2D}
}%
\hfill
\subfigure[Elliptic orbit of $m_2$ relative to $m_1$]{%
\centering
    \includegraphics[width=0.6\columnwidth, trim={0.1in 0in 0.1in 0in},clip]{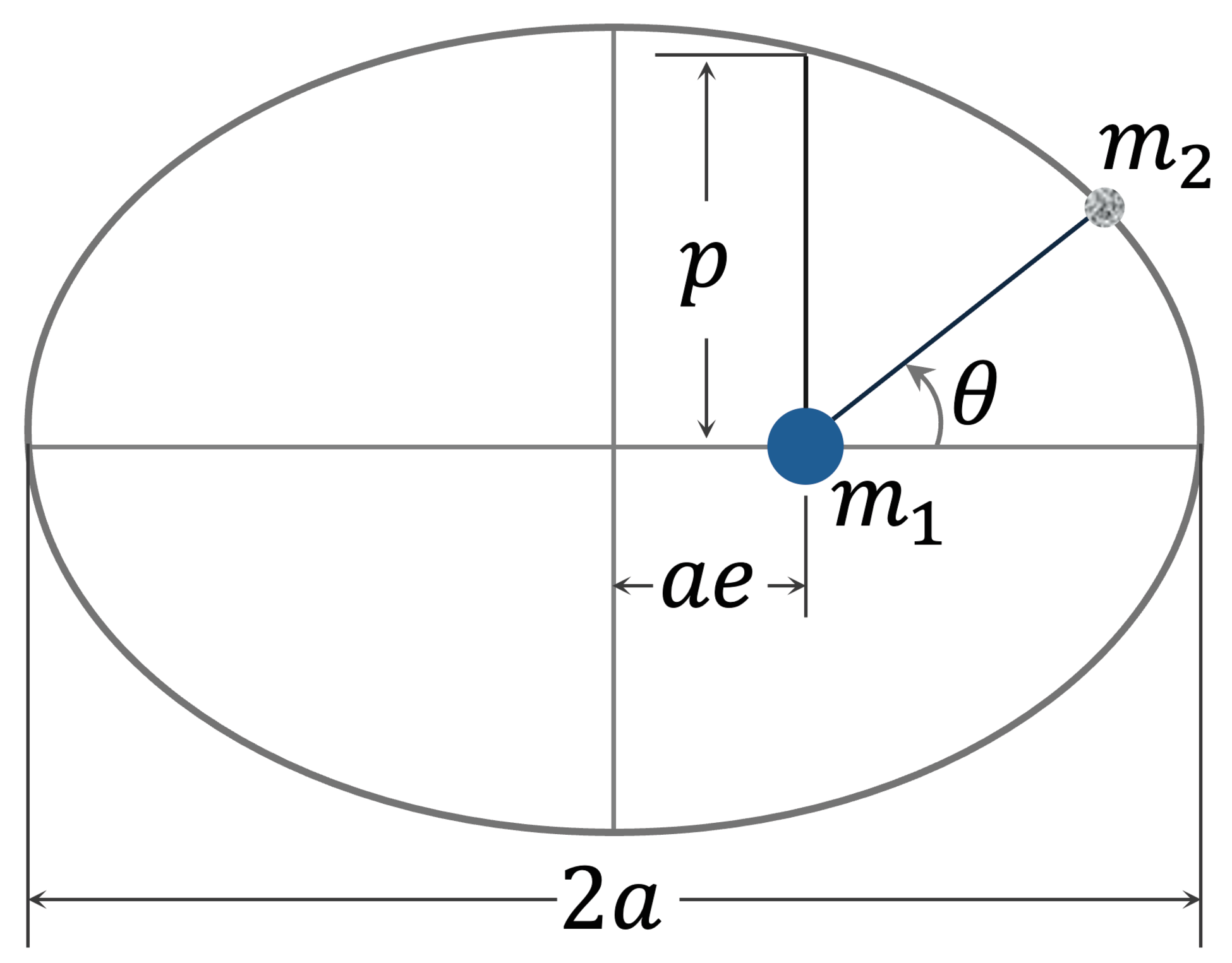}
    \label{fig:Orbit2Rel1}
}%
\caption{ER3BP coordinate frames and various elliptic orbits}
\label{fig: ColocatedFrames}
\end{figure}

Figure~\ref{fig: ColocatedFrames} depicts the two primaries that are assumed to be rotating about their common center of mass in an elliptic orbit.
%
%
The unit vectors $\hat{\bi}$, $\hat{\bj}$, $\hat{\bk}$ are the orthogonal directions in the local frame $L$ which rotates around its $\bar{z}$ axis with respect to the inertial frame $N$ at a non-uniform angular velocity of $\bomega_{L/N}$,
\begin{equation}\label{eq:omega=thetadot}
\bomega_{L/N} = \omega \hat{\bk} = \dot\theta \hat{\bk}
\end{equation}
where $\theta$ is the true anomaly.  In this paper, we forgo some standard notation in orbital mechanics such as $f$ for true anomaly in favor of saving such symbols for other prevalent usage such as ``$f$ for function.''   The distance between the two primaries, $\xi$, varies as a function of the true anomaly, $\theta$,
\begin{equation}
    \xi = \frac{a(1-e^2)}{1 + e\cos{\theta}} = \frac{p}{1 + e\cos{\theta}}
    \label{eq: ERTBP_ellipticTimeVaryingDistance}
\end{equation}
where $a$ is the relative semi-major axis (see Fig.~\ref{fig:Orbit2Rel1}), $e$ is the eccentricity and $p$ is the semilatus rectum of the relative elliptic orbit.  That is, $a$ and $p$ correspond to the primaries' relative orbit about each other as shown in Fig.~\ref{fig:Orbit2Rel1}, and not the elliptic orbits shown in Fig.~\ref{fig: CoordFrame_BaryToInertial_2D}. Not shown in Fig.~\ref{fig: ColocatedFrames} is the relative elliptic orbit of $m_1$ with respect to $m_2$. The orbits depicted in Figs.~\ref{fig: CoordFrame_BaryToInertial_2D} and \ref{fig:Orbit2Rel1} are geometrically similar such that the barycentric semi-major axes $a_1$ and $a_2$ of  $m_1$ and $m_2$ respectively (see Fig.~\ref{fig: CoordFrame_BaryToInertial_2D}) are given by\cite{szebehely_theory_1967},
$$a_1 = a \mu \quad a_2 = a(1-\mu)  $$
The eccentricities of all of the orbits shown in Figs.~\ref{fig: CoordFrame_BaryToInertial_2D} and \ref{fig:Orbit2Rel1} including the relative elliptic orbit of $m_1$ with respect to $m_2$ are the same\cite{szebehely_theory_1967}.

\subsection{Basic Equations of Motion}
The primaries are located along the $x$ axis at locations $(\bar{x}_1,\,0,\,0)$ and $(\bar{x}_2,\,0,\,0)$. We use the ``overbar'' notation to express quantities in the local frame $L$ in order to save the unbarred notation for the Nechvile-transformed variables defined later in this section. The vector $\bar\br = (\bar x, \bar y, \bar z)$ is the position vector of the spacecraft with respect to the barycenter expressed in $L$. The position of the spacecraft in the local, rotating frame relative to the primaries are given by the vectors $\bar{\br}_1$ and $\bar{\br}_2$. The magnitudes of these vectors are given by,
\begin{equation}\label{eq: rho1Andrho2}
\begin{aligned}
 \norm{\bar\br_1}_2 =   \bar{\rho}_1 := \sqrt{(\bar{x}-\bar{x}_1)^2 + \bar{y}^2 + \bar{z}^2} \\ \norm{\bar\br_2}_2 = \bar{\rho}_2 := \sqrt{(\bar{x}-\bar{x}_2)^2 + \bar{y}^2 + \bar{z}^2}
\end{aligned}
\end{equation}
Consistent with these assumptions\cite{elipe-wrongHyp-2024}, gravitational forces are modeled using Newton's inverse-square law.  A generic external force $\bF$ acts only on the spacecraft (and not on the two primaries). The components of $\bF$ in $L$ are given by $F_x$, $F_y$ and $F_z$.

By a straightforward application of Newton's Laws of motion and gravity we get (see Fig.~\ref{fig: ColocatedFrames}),
\begin{equation}\label{eq:F=ma}
m_{sc}\frac{d^2\bar{\br}}{dt^2}\bigg|_N = -G\frac{m_{sc}m_1}{\norm{\bar\br_1}_2^3}\bar\br_1  -G\frac{m_{sc}m_2}{\norm{\bar\br_2}_2^3}\bar\br_2 + \bF
\end{equation}
where the subscript $N$ in \eqref{eq:F=ma} denotes derivatives with respect to the inertial frame.  From basic kinematics, we can express the second time derivative of $\bar\br$ as,
\begin{equation}
    \begin{aligned}
        \frac{d^2\bar{\br}}{dt^2}\bigg|_N &= \frac{d^2\bar{\br}}{dt^2}\bigg|_L + \frac{d\bomega}{dt} \times \bar{\br} + 2\bomega \times \frac{d\bar{\br}}{dt}\bigg|_L + \bomega \times \left(\bomega \times \bar{\br} \right)
    \end{aligned}
    \label{eq: ERTBP_2nd_inertialDerivative}
\end{equation}
Substituting \eqref{eq:omega=thetadot} in \eqref{eq: ERTBP_2nd_inertialDerivative}, we get,
\begin{equation}
    \ddot{\bar{\br}}\big|_N = \ddot{\bar{x}}\,\hat{\bi} + \ddot{\bar{y}}\hat{\bj} + \ddot{\bar{z}}\,\hat{\bk} - \dot{\omega}\,\bar{y}\,\hat{\bi} + \dot{\omega}\,\bar{x}\hat{\bj} - 2\omega\,\dot{\bar{y}}\,\hat{\bi} + 2\omega\,\dot{\bar{x}}\hat{\bj} - \omega^2\,\bar{x}\,\hat{\bi} - \omega^2\,\bar{y}\hat{\bj}
    \label{eq: ERTBP_2nd_inertialDerivative_ScalarForm}
\end{equation}
Substituting \eqref{eq: ERTBP_2nd_inertialDerivative_ScalarForm} in \eqref{eq:F=ma} and equating components we arrive at,
\begin{subequations}\label{eq: ERTBP_forceBalance_dimensional}
    \begin{gather}
\ddot{\bar{x}} - \dot{\omega}\,\bar{y} - 2\,\omega\,\dot{\bar{y}} - \omega^2 \bar{x} = -G\,m_1\frac{\bar{x}-{\bar{x}_1}}{\bar{\rho}_1^3} - G\,m_2\frac{\bar{x}-{\bar{x}_2}}{\bar{\rho}_2^3} + \frac{F_{x}}{m_{sc}} \label{eq: ERTBP_XforceBalance_dimensional}
\\
\ddot{\bar{y}} + \dot{\omega}\,\bar{x} + 2\,\omega\,\dot{\bar{x}} - \omega^2 \bar{y} = -G\,m_1\frac{\bar{y}}{\bar{\rho}_1^3} - G\,m_2\frac{\bar{y}}{\bar{\rho}_2^3} + \frac{F_{y}}{m_{sc}} \label{eq: ERTBP_YforceBalance_dimensional}
\\
\ddot{\bar{z}} = -G\,m_1\frac{\bar{z}}{\bar{\rho}_1^3} - G\,m_2\frac{\bar{z}}{\bar{\rho}_2^3} + \frac{F_{z}}{m_{sc}} \label{eq: ERTBP_ZforceBalance_dimensional}
\end{gather}
\end{subequations}

\subsection{Nechvile’s Transformation}

Nechvile's transformation\cite{szebehely_theory_1967,broucke-nechville-1966,duboshin-1971} involves scaling the components of the spacecraft position vector by \eqref{eq: ERTBP_ellipticTimeVaryingDistance} and changing the independent variable from time to true anomaly in \eqref{eq: ERTBP_forceBalance_dimensional}.  This results in
\begin{equation}
    x := \frac{\bar{x}}{\xi}, \quad y := \frac{\bar{y}}{\xi}, \quad z := \frac{\bar{z}}{\xi}
    \label{eq: ERTBP_rotatingFrame_distanceScaling}
\end{equation}
and
\begin{equation}
    \frac{d}{dt}=\omega\,\frac{d}{d\theta}, \qquad \frac{d^2}{dt^2} = \dot{\omega}\, \frac{d}{d\theta} + \omega^2\, \frac{d^2}{d\theta^2}
    \label{eq: ERTBP_trueAnomalyChainRule}
\end{equation}
An immediate consequence of \eqref{eq: ERTBP_rotatingFrame_distanceScaling} is that the rotating frame is also pulsating. In other words, Nechvile's frame is rotating and pulsating. In this frame,  the position of the primaries are fixed and given by $(-\mu,\,0,\,0)$ and $(1-\mu,\,0,\,0)$.  Conversely, in the non-pulsating rotating frame, the primaries pulsate nonuniformly with orbital frequency $\omega$.

Using \eqref{eq: ERTBP_rotatingFrame_distanceScaling} and \eqref{eq: ERTBP_trueAnomalyChainRule}, the left hand side of \eqref{eq: ERTBP_forceBalance_dimensional} transforms to,
\begin{subequations}\label{eq: ERTBP_LHS_transformed}
\begin{gather}
    2\omega\,x'\left( \frac{d\xi}{dt}\right) + \xi\,\dot{\omega}\,x' + \xi\,\omega^2\,x'' + x\left(\frac{d^2\xi}{dt^2}\right)
    -\dot{\omega}\,y\,\xi-2\omega^2\xi\,y' \nonumber \\
    - 2\omega\,y\left(\frac{d\xi}{dt}\right) - \omega^2\xi\,x
    \label{eq: ERTBP_xDirection_LHS_transformed}
    \\
    2\omega\,y'\left( \frac{d\xi}{dt}\right) + \xi\,\dot{\omega}\,y' + \xi\,\omega^2\,y'' + y\left(\frac{d^2\xi}{dt^2}\right)+\dot{\omega}\,x\,\xi + 2\omega^2\xi\,x' \nonumber \\
    + 2\omega\,x\left(\frac{d\xi}{dt}\right) - \omega^2\xi\,y
    \label{eq: ERTBP_yDirection_LHS_transformed}
    \\
    2\omega\,z'\left(\frac{d\xi}{dt}\right) + \xi\,\dot{\omega}\,z' + \xi\,\omega^2\,z'' + z\left(\frac{d^2\xi}{dt^2}\right)
    \label{eq: ERTBP_zDirection_LHS_transformed}
\end{gather}
\end{subequations}
The equation of motion between the two primaries,
\begin{equation}
    \frac{d^2\xi}{dt^2} -\xi\omega^2 = -G\frac{m_1 + m_2}{\xi^2}
    \label{eq: ERTBP_twoBodyEOM}
\end{equation}
and the conservation of angular momentum in the two-body problem,
\begin{equation}
\xi\dot{\omega} + 2\frac{d\xi}{dt}\,\omega = 0
    \label{eq: ERTBP_conservationAngMom_Derivative}
\end{equation}
facilitate an elimination of the first- and second time-derivatives of $\xi$ in \eqref{eq: ERTBP_LHS_transformed}.  After \eqref{eq: ERTBP_twoBodyEOM} and \eqref{eq: ERTBP_conservationAngMom_Derivative} are used to carry out this elimination step, \eqref{eq: ERTBP_LHS_transformed}, and hence the left-hand-side of \eqref{eq: ERTBP_forceBalance_dimensional}, simplifies to,
\begin{subequations}\label{eq: ERTBP_LHS_transform_rewrite2nd}
\begin{gather}
    \xi\,\omega^2\left(x'' - 2y'\right) - x \frac{\xi^2\omega^2}{p}
    \label{eq: ERTBP_xDirection_LHS_transform_rewrite2nd}
    \\
    \xi\,\omega^2(y'' + 2x') - y\frac{\xi^2\omega^2}{p}
    \label{eq: ERTBP_yDirection_LHS_transform_rewrite2nd}\\
    \xi\,\omega^2 \, z'' + z \left[ \xi\,\omega^2 - \frac{\xi^2\omega^2}{p} \right]
    \label{eq: ERTBP_zDirection_LHS_transform_rewrite2nd}
\end{gather}
\end{subequations}
In the derivation of \eqref{eq: ERTBP_LHS_transform_rewrite2nd}, we have also used the connection between the semilatus rectum $p$ and the angular momentum $h$:
\begin{equation}\label{eq: ERTBP_conservationAngularMomentum}
p = \frac{h^2}{G(m_1 + m_2)} = \frac{(\xi^2 \omega)^2}{G(m_1 + m_2)}
\end{equation}

\subsection{The Transformed Equations of Motion}
Substituting \eqref{eq: ERTBP_LHS_transform_rewrite2nd} for the left-hand-side of \eqref{eq: ERTBP_forceBalance_dimensional} results in,
\begin{subequations}\label{eq: ERTBP_LHSandRHS_withScaleFactor}
\begin{align}
    \xi\,\omega^2\left(x'' - 2y'\right) - x \frac{\xi^2\omega^2}{p} &= -G\,m_1\frac{x-{x}_1}{\xi^2\rho_1^3} \nonumber \\
    &\quad - G\,m_2\frac{x-{x}_2}{\xi^2\rho_2^3} + \frac{F_x}{m_{sc}}
    \label{eq: ERTBP_xDirection_LHSandRHS_withScaleFactor}
    \\
    \xi\,\omega^2(y'' + 2x') - y\frac{\xi^2\omega^2}{p} &= -G\,m_1\frac{y}{\xi^2\rho_1^3} \nonumber \\
   &\quad - G\,m_2\frac{y}{\xi^2\rho_2^3} + \frac{F_y}{m_{sc}}
    \label{eq: ERTBP_yDirection_LHSandRHS_withScaleFactor}
    \\
    \xi\,\omega^2z'' + z\left[ \xi\,\omega^2 - \frac{\xi^2\omega^2}{p} \right] &= -G\,m_1\frac{z}{\xi^2\rho_1^3} \nonumber \\
    &\quad - G\,m_2\frac{z}{\xi^2\rho_2^3} + \frac{F_z}{m_{sc}}
    \label{eq: ERTBP_zDirection_LHSandRHS_withScaleFactor}
\end{align}
\end{subequations}
where
$$\rho_1 := \frac{\bar\rho_1}{\xi}, \quad  \rho_2 := \frac{\bar\rho_2}{\xi} $$
Using \eqref{eq: massParameter_mu} and \eqref{eq: ERTBP_conservationAngularMomentum},
\eqref{eq: ERTBP_LHSandRHS_withScaleFactor} can be rewritten as,
\begin{subequations}\label{eq: ERTBP_LHSandRHS_conAngMomVersion}
    \begin{align}
    x'' - 2y' &= \frac{\xi}{p}\left[x - \frac{(1-\mu)(x-x_1)}{\rho_1^3} - \frac{\mu\,(x-x_2)}{\rho_2^3} \right]+ \frac{F_x}{m_{sc}\,\omega^2\xi}
    \label{eq: ERTBP_xDirection_LHSandRHS_conAngMomVersion}
    \\
    y'' + 2x' &= \frac{\xi}{p} \left[ y - \frac{(1-\mu)\,y}{\rho_1^3} - \frac{\mu\,y}{\rho_2^3} \right]+ \frac{F_y}{m_{sc}\,\omega^2\xi}
    \label{eq: ERTBP_yDirection_LHSandRHS_conAngMomVersion}
    \\
    z'' &= -z + \frac{\xi}{p} \left[ z - \frac{(1-\mu)\,z}{\rho_1^3} - \frac{\mu\,z}{\rho_2^3} \right]+ \frac{F_z}{m_{sc}\,\omega^2\xi}
    \label{eq: ERTBP_zDirection_LHSandRHS_conAngMomVersion}
\end{align}
\end{subequations}
Substituting \eqref{eq: ERTBP_ellipticTimeVaryingDistance} in \eqref{eq: ERTBP_conservationAngularMomentum}, the centrifugal term $\omega^2\,\xi$ can be rewritten as:
\begin{equation}
    \omega^2\,\xi = \frac{p\,G\,(m_1+m_2)}{\xi^3} = \frac{G\,(m_1+m_2)}{p^2}\,(1+e\cos{\theta})^3
    \label{eq: ERTBP_omegaSqXi_scale}
\end{equation}
Combining \eqref{eq: ERTBP_LHSandRHS_conAngMomVersion} and \eqref{eq: ERTBP_omegaSqXi_scale} generates the following equations:
\begin{subequations}\label{eq: ERTBP_xyzDirection_LHSandRHS_SzebWithThrust}
\begin{align}
    x'' - 2y' &= \frac{1}{1+e\cos{\theta}}\left[x - \frac{(1-\mu)(x-x_1)}{\rho_1^3} - \frac{\mu\,(x-x_2)}{\rho_2^3} \right]  \nonumber \\
    &\quad + \frac{p^2}{(1+e\cos{\theta})^3} \frac{F_x/m_{sc}}{G\,(m_1+m_2)}
    \\
    y'' + 2x' &= \frac{1}{1+e\cos{\theta}} \left[ y - \frac{(1-\mu)\,y}{\rho_1^3} - \frac{\mu\,y}{\rho_2^3} \right] \nonumber \\
   &\quad +
    \frac{p^2}{(1+e\cos{\theta})^3} \frac{F_y/m_{sc}}{G\,(m_1+m_2)}
    \\
    z'' &= \frac{1}{1+e\cos{\theta}} \left[ -ze\cos{\theta} - \frac{(1-\mu)\,z}{\rho_1^3} - \frac{\mu\,z}{\rho_2^3} \right] \nonumber \\
    &\quad + \frac{p^2}{(1+e\cos{\theta})^3} \frac{F_z/m_{sc}}{G\,(m_1+m_2)}
    \end{align}
\end{subequations}
Let $\Psi$ be the pseudo-potential function that governs the gravitational and centrifugal forces defined by\cite{szebehely_theory_1967}:
\begin{multline}
    \Psi(x,y,z, \theta) := \frac{1}{1+e\cos{\theta}}\left[\frac{1}{2}(x^2+y^2-z^2e\cos{\theta}) \right. \\ 
    \left. + \frac{1-\mu}{\rho_1(x,y,z)} + \frac{\mu}{\rho_2(x,y,z)} \right]
    \label{eq: ERTBP_PseudoPotentialFunction}
\end{multline}
Then \eqref{eq: ERTBP_xyzDirection_LHSandRHS_SzebWithThrust} can be written  more succinctly as,
\begin{subequations}
    \begin{align}
         x'' - 2y' &= \partial_{x}\Psi(x,y,z, \theta) + \frac{1}{(1+e\cos{\theta})^3}\,   \frac{F_x/m_{sc}}{A_C}
         \\
         y'' + 2x' &= \partial_{y}\Psi(x,y,z, \theta) + \frac{1}{(1+e\cos{\theta})^3}\,\frac{F_y/m_{sc}}{A_C}
         \\
         z'' &= \partial_{z}\Psi(x,y,z, \theta) + \frac{1}{(1+e\cos{\theta})^3}\,\frac{F_z/m_{sc}}{A_C}
    \end{align}
    \label{eq: ERTBP_EOM_pseudoPot_ver}
\end{subequations}
where $A_C$ is a constant defined by,
\begin{equation}
    A_C:= \frac{G\,(m_1+m_2)}{p^2}
    \label{eq: ERTBP_accelScaleFactor}
\end{equation}
Equation~(\ref{eq: ERTBP_EOM_pseudoPot_ver}) (or \eqref{eq: ERTBP_xyzDirection_LHSandRHS_SzebWithThrust}) with $F_x = 0= F_y = F_z$ are well known\cite{szebehely_theory_1967,broucke_stability_1969,celletti_dynamics_2024,peng-2015}.  Adding an external force to these equations generates the cubic multiplicative term introduced in \eqref{eq:cubic}. For $e = 0$, this term is unity. However, when $e \ne 0$, the cubic term has far reaching consequences.

\subsection{Constant and Nechvile Units}
\begin{figure}[h!]
\centering
\subfigure[A 2:1 Halo orbit in the Nechvile (rotating and pulsating) frame]{%
    \includegraphics[width=0.45\textwidth, trim={1.3in 0 1in 0in},clip]{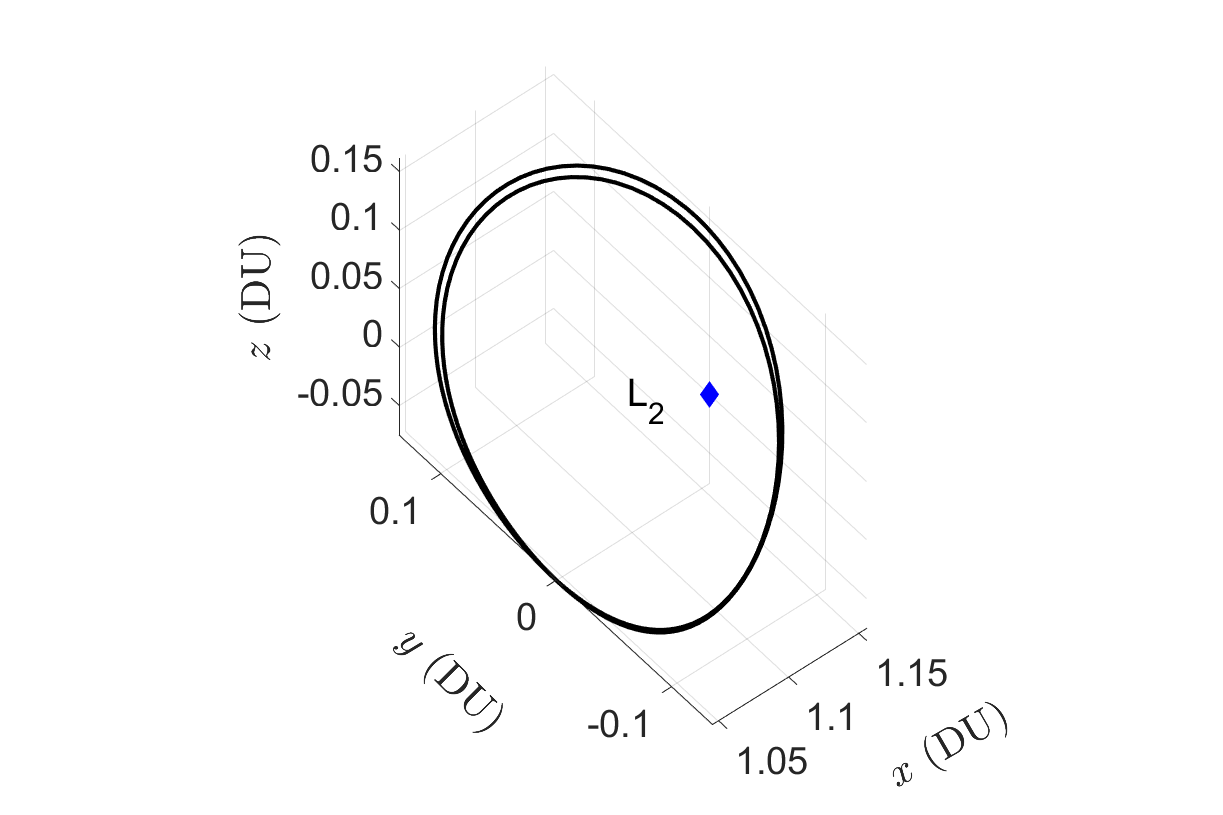}
    \label{fig:Halo-a}
}%
\hfill
\subfigure[Halo orbit of Fig.~\ref{fig:Halo-a} in a rotating only (nonpulsating) frame ]{%
    \includegraphics[width=0.45\textwidth, trim={1.3in 0 1in 0in},clip]{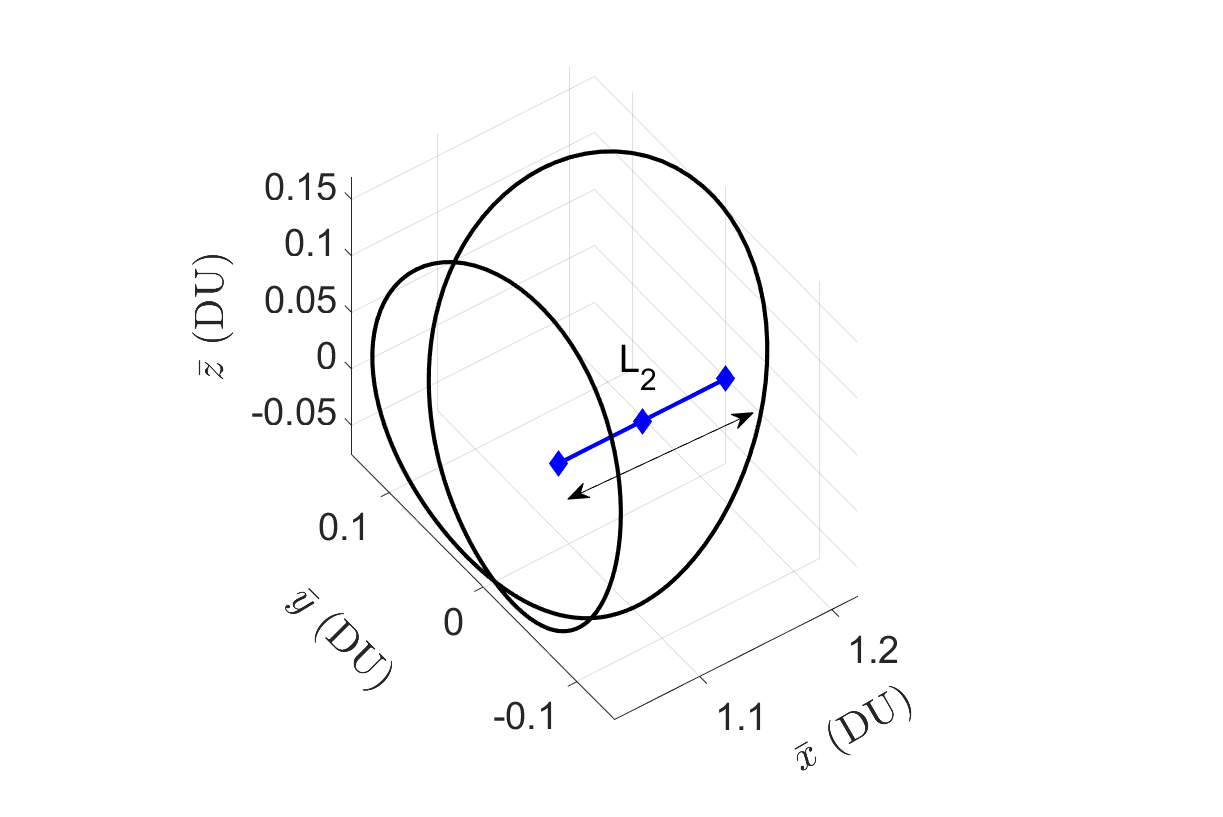}
    \label{fig:Halo-b}
}%
\caption{Halo orbits and $L_2$ points in pulsating and nonpulsating frames}
\label{fig:Halos}
\end{figure}

Based on Eqs~(\ref{eq: ERTBP_ellipticTimeVaryingDistance}), (\ref{eq: ERTBP_rotatingFrame_distanceScaling}) and (\ref{eq: ERTBP_accelScaleFactor}), 
we define a constant distant unit, $D_C$, to be the same as the semilatus rectum,
\begin{equation}\label{eq:D0=bydef}
D_C := p
\end{equation}
and a Nechvile distant unit $D_N$ ($= \xi$) to be defined by,
\begin{equation}
    D_N := \frac{D_C}{1+e\cos{\theta}}
    \label{eq: ERTBP_DU_define}
\end{equation}
Obviously, the Nechvile distant unit is time-varying. From \eqref{eq: ERTBP_omegaSqXi_scale}, we have,
\begin{equation}
     \omega = (1+e\cos{\theta})^2 \, \sqrt{\frac{G\,(m_1+m_2)}{p^3}}
    \label{eq: ERTBP_OM_define}
\end{equation}
The square-root term in \eqref{eq: ERTBP_OM_define} is the constant angular velocity of a fictitious circular orbit of radius $D_C$; hence, we define,
\begin{equation}\label{eq:omega=funtheta}
     \omega_C:= \sqrt{\frac{G\,(m_1+m_2)}{p^3}}, \qquad \omega_N := (1+e\cos{\theta})^2 \, \omega_C
\end{equation}
and refer to $\omega_N$ ($= \omega$) as the (time-varying) Nechvile frequency.  The Nechvile unit of time $\tau_N$ is therefore ``derived'' by,
\begin{equation}
    \tau_N := \frac{1}{\omega_N}  = \frac{\tau_C}{(1+e\cos{\theta})^2}; \quad \tau_C := \frac{1}{\omega_C}
    \label{eq: ERTBP_TU_define}
\end{equation}
Note that the Nechvile unit of time as defined in \eqref{eq: ERTBP_TU_define} is also consistent with the change in the unit of time as implied by \eqref{eq: ERTBP_trueAnomalyChainRule}.  We briefly note that the motion of a spacecraft in Nechvile's pulsating coordinates is not the same as its ``physical motion'' in non-pulsating coordinates. For instance, a 2:1 Halo orbit around the Earth-Moon $L_2$ point is shown in Fig.~\ref{fig:Halo-a} in Nechvile's coordinates. The same orbit in a non-pulsating (but rotating) coordinate system is shown in Fig~\ref{fig:Halo-b}. Note also that the Lagrange point ($L_2$) is not fixed in the (non-pulsating) rotating frame. As shown in Fig.~\ref{fig:Halo-b} the $L_2$ point oscillates at the orbital frequency. These interesting features of the ER3BP are in sharp contrast to motion in the CR3BP\cite{peng-2015}.
%

Using Eqs.~(\ref{eq:D0=bydef}), (\ref{eq: ERTBP_DU_define}) and (\ref{eq: ERTBP_TU_define}), we can derive constant and Nechvile units of acceleration as,
\begin{equation}\label{eq:AN=bydef}
A_C := \frac{D_C}{\tau_C^2}, \qquad A_N := \frac{D_N}{\tau_N^2} = A_C (1+e\cos{\theta})^3
\end{equation}
The quantity $A_C$ defined in \eqref{eq:AN=bydef} is exactly the same as that defined in \eqref{eq: ERTBP_accelScaleFactor} with the former being derived ``independently'' by consideration of a fictitious circular orbit.  From similar considerations, we define,
\begin{equation}\label{eq:vC=}
v_C := \frac{A_C}{\omega_C} = \sqrt{\frac{G(m_1+m_2)}{p}}
\end{equation}

\section{ER3BP Models For a Thrusting Spacecraft}
\label{sec:discuss}

The rocket mass flow rate equation is given by,
\begin{equation}
    \dot{m}_{sc} = -\frac{\|\bT\|_p}{I_{sp}\,g_0}
    \label{eq: ERTBP_timeMassFlowRate}
\end{equation}
where $\|\bT\|_p$ is the $\ell^p$-norm of $\bT \in \Real^3$, $p \in \set{1, 2, \infty} $\cite{ross_how_2004,ross_space_2006}, $I_{sp}$ is the specific impulse of the propellant and $g_0$ is the gravitational acceleration at Earth sea level.  Changing the independent variable in \eqref{eq: ERTBP_timeMassFlowRate} from time to true anomaly and using Eqs.~(\ref{eq:omega=thetadot})  and (\ref{eq: ERTBP_OM_define})  we get,
\begin{equation}
     m_{sc}' = -\frac{1}{\omega}\,\frac{\|\bT\|_p}{I_{sp}\,g_0}
    \label{eq: ERTBP_trueAnomMassFlowRate}
\end{equation}
Substituting \eqref{eq: ERTBP_OM_define} and \eqref{eq:omega=funtheta} in \eqref{eq: ERTBP_trueAnomMassFlowRate} yields,
\begin{equation}
     m_{sc}' = -\frac{1}{(1+e\cos{\theta})^2}\frac{\|\bT\|_p}{I_{sp}\,\omega_C \,g_0}
    \label{eq: ERTBP_trueAnomMassFlowRate_withCharTime}
\end{equation}
Let $m_C$ be a constant mass associated with the spacecraft such as its initial value. Define a dimensionless spacecraft mass variable as,
\begin{equation}
    m:=\frac{ m_{sc}}{m_C}
    \label{eq: ERTBP_nonDmass}
\end{equation}
Substituting \eqref{eq: ERTBP_nonDmass} in \eqref{eq: ERTBP_trueAnomMassFlowRate_withCharTime} and augmenting the resulting equation with \eqref{eq: ERTBP_EOM_pseudoPot_ver}, we get the following model (hereafter referred to as Model $A$) for the basic equations of motion for a thrusting spacecraft with time-varying mass:
\begin{subequations} \label{eq: ERTBP_EOM_pseudoPot_ver_withNonDmass}
\begin{align}
         x'' - 2y' &= \partial_{x}\Psi(x,y,z, \theta) + \frac{1}{(1+e\cos{\theta})^3}\,   \frac{T_x/ (m_C\,A_C)}{ m}
         \\
         y'' + 2x' &= \partial_{y}\Psi(x,y,z, \theta) + \frac{1}{(1+e\cos{\theta})^3}\,\frac{T_y/ (m_C\,A_C)}{m }
         \\
         z'' &= \partial_{z}\Psi(x,y,z, \theta) + \frac{1}{(1+e\cos{\theta})^3}\,\frac{T_z/ (m_C\,A_C)}{m }\\
         m' &= -\frac{1}{(1+e\cos{\theta})^2}\frac{\|\bT\|_p}{({I}_{sp}\,\omega_C)(m_C\,g_0)}
\end{align}
\end{subequations}
Note the ``discrepancy'' between the quadratic term in the rocket mass flow rate equation and the cubic terms in remainder of \eqref{eq: ERTBP_EOM_pseudoPot_ver_withNonDmass}. As explained in Section~\ref{sec:intro}, the consequences of this discrepancy are far reaching.

\subsection{Models $\bA$ Versus $\bB$}
Suppose we scale the thrust force by the cubic term in \eqref{eq: ERTBP_EOM_pseudoPot_ver_withNonDmass} in order to ``simplify'' Model $A$; then, \eqref{eq: ERTBP_EOM_pseudoPot_ver_withNonDmass} simplifies to the following Model $B$:
\begin{subequations} \label{eq: ERTBP_EOM_pseudoPot_ver_withNechville}
\begin{align}
         x'' - 2y' &= \partial_{x}\Psi(x,y,z, \theta) +    \frac{\widetilde{T}_x/(m_C\,A_C)}{m}
         \\
         y'' + 2x' &= \partial_{y}\Psi(x,y,z, \theta) + \frac{\widetilde{T}_y/(m_C\,A_C)}{m}
         \\
         z'' &= \partial_{z}\Psi(x,y,z, \theta) + \frac{\widetilde{T}_z/(m_C\,A_C)}{m}\\
         m' &= -{(1+e\cos{\theta})}\frac{\|\widetilde{\bT}\|_p}{({I}_{sp}\,\omega_C)(m_C\,g_0)}
\end{align}
\end{subequations}
where $ \widetilde{\bT} $ (and similarly, $\widetilde{T}_x$, $\widetilde{T}_y$ and $\widetilde{T}_z$) are defined by,
\begin{equation}\label{eq:Ttilde=bydef}
\widetilde{\bT} := \frac{\bT}{(1+e\cos{\theta})^3}
\end{equation}
Because Model $B$ is equivalent to scaling the thrust force using Nechvile's acceleration units, it has the appearance of simplicity.  Note, however, that although the cubic term in Model $A$ has disappeared in Model $B$, the mass-flow equation in the latter now contains the linear cosine term, $(1 + e \cos\theta)$.

Obviously, Models $A$ and $B$ are mathematically and physically equivalent; however,  there are practical consequences and differences in using \eqref{eq: ERTBP_EOM_pseudoPot_ver_withNonDmass} over \eqref{eq: ERTBP_EOM_pseudoPot_ver_withNechville} in terms of the controlled motion of a spacecraft, particularly for trajectory optimization and guidance.  This is because in a practical spacecraft\cite{biblarz}, the magnitude of the thrust force is bounded by some finite value, $T_{max}$. This implies\cite{ross_how_2004,ross_space_2006,ross-book},
\begin{equation}\label{eq:Tq-bound}
\norm{\bT}_q \le T_{max}, \qquad \frac{1}{q} + \frac{1}{p} = 1
\end{equation}
To limit the scope of the discussion to follow, we consider the simple symmetric case of $p=q= 2$ corresponding to a single thruster with thrust force $T \ge 0$. Then,
\begin{subequations}\label{eq:TboundsA+B}
\begin{align}
T  & \le T_{max} \label{eq:TboundA}\\
\Rightarrow \quad \widetilde{T} &\le \frac{T_{max}}{(1+e\cos{\theta})^3} \label{eq:TboundB}
\end{align}
\end{subequations}
From \eqref{eq:TboundB}, it follows that the control space (i.e., allowable values of the control\cite{ross-book}) for Model $B$ is time-varying (with $\widetilde{T}$ as a control variable) while that for Model $A$ (Cf.~\eqref{eq:TboundA}) is time-invariant (i.e., with $T$ as a control variable).
In the preceding sentence and the discussions to follow we use ``time'' interchangeably with ``true anomaly'' because $\theta$ is the independent variable.
Thus, the purported simplicity of the dynamical equations offered by Model $B$ via \eqref{eq: ERTBP_EOM_pseudoPot_ver_withNechville} comes at the price of intricacies in its time-varying control space governed by \eqref{eq:TboundB}. We hasten to note that Pontryagin's Principle in optimal control theory can be just as easily applied to time-varying control spaces as time-invariant ones\cite{ross-book}; however, any computational advantage offered by the ``simplicity'' of the dynamics of Model~$B$ is lost in its transformation of the control space from time-invariant to time-varying.  In other words, the ``problem'' has not been solved; rather, it has simply been shifted from the dynamical equations to time-variability in the control space.

\subsection{Model $\widetilde{\bB}$ as an Alternative to Model $\bB$}

\begin{figure}[h!]
\centering
\subfigure[Reduced thruster capacity for Model~$\widetilde{B}$ for $e=0.2$]{%
    \includegraphics[width=0.45\textwidth]{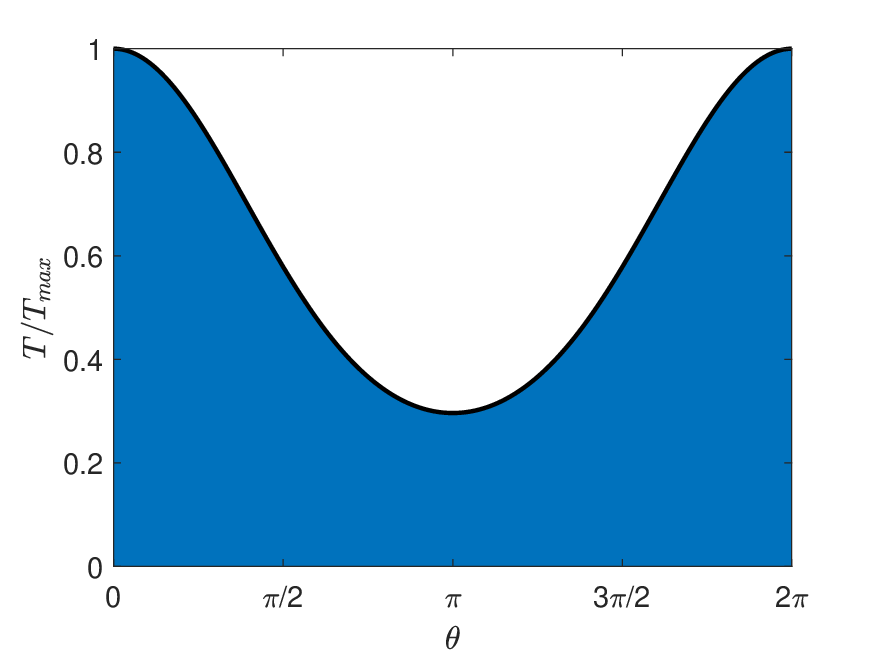}
    \label{fig:ModelCspace}
}%
\hfill
\subfigure[Worst-case limitation of thruster capacity for Model~$\widetilde{B}$ ]{%
    \includegraphics[width=0.45\textwidth]{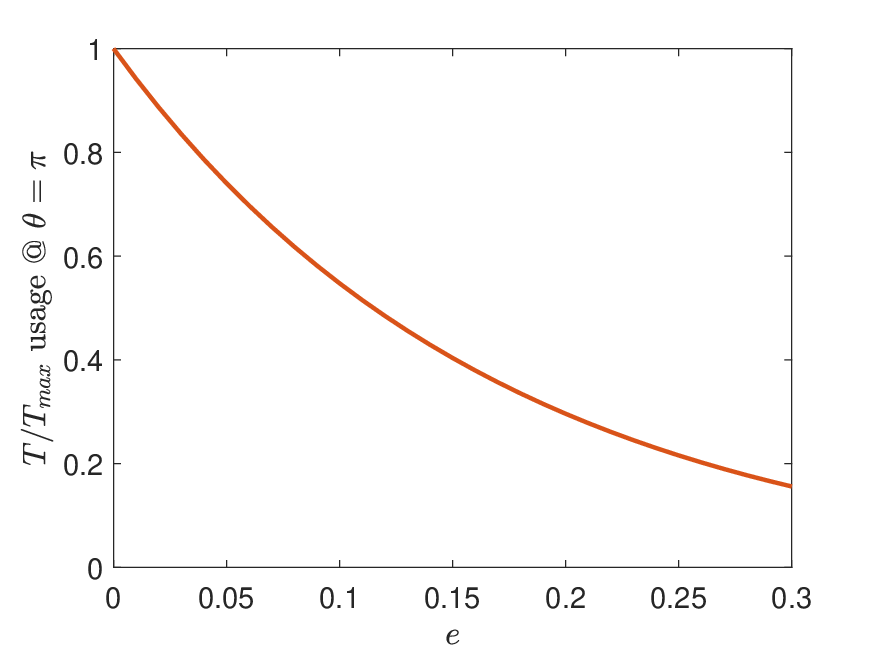}
    \label{fig:ModelCworstLimits}
}%
\caption{Limitations on maximum thrust usage imposed by Model~$\widetilde{B}$  }
\label{fig:ModelCplots}
\end{figure}

To overcome the variability of the control space given by \eqref{eq:TboundB}, suppose we simplified Model $B$ to construct Model~$\widetilde{B}$ wherein the dynamics is the same as that of Model $B$ (i.e., \eqref{eq: ERTBP_EOM_pseudoPot_ver_withNechville}) but the control space is set to be time-invariant as in Model~$A$.  That is, in Model $\widetilde{B}$ we simply set,
\begin{equation}\label{eq:ModelC}
\widetilde{T }  \le \widetilde{T}_{max}
\end{equation}
where $\widetilde{T}_{max}$ is a constant. If $\widetilde{T}_{max}$ is set to $T_{max}$, it is quite possible that a feasible solution for the actual thrust $T$ might not exist.  This is because according to \eqref{eq:Ttilde=bydef}, a constant value of $\widetilde{T}_{max}$ maps to a modulated value of $T$ according to,
\begin{equation}\label{eq:Ttilde2T}
T = \widetilde{T}_{max} (1+e\cos{\theta})^3
\end{equation}
Hence, at $\theta = 0$, $T$ will exceed $T_{max}$ if $\widetilde{T}_{max} = T_{max}$.  In fact, the spacecraft will not be able to deliver  $\widetilde{T}_{max} = T_{max}$ whenever $(1+e\cos{\theta})^3 > 1$. It is a simple matter to show that this condition is invariant of $e$ and can be expressed as,
$$\theta \in \set{[0, \pi/2)\cup (3\pi/2, 1]  }$$
Thus, the spacecraft will not be able to deliver $\widetilde{T}_{max} = T_{max}$
over half an orbit. In order to meet the physical requirement of $T\le T_{max}$, we must have,
\begin{equation}\label{eq:maxTtildeCondition-pre}
 \max_{\theta} \left\{ \widetilde{T}_{max} (1+e\cos{\theta})^3 \right\}  \le T_{max}
\end{equation}
Equation~(\ref{eq:maxTtildeCondition-pre}) simplifies to,
\begin{equation}\label{eq:maxTtildeCondition}
\widetilde{T}_{max}  \le \frac{T_{max}}{(1+e)^3}
\end{equation}
The implications of \eqref{eq:maxTtildeCondition} are quite severe.  Mapping \eqref{eq:maxTtildeCondition} back to ``$T$-space,'' we get (from Eqs.~(\ref{eq:Ttilde=bydef}), (\ref{eq:ModelC}) and (\ref{eq:maxTtildeCondition})),
\begin{equation}\label{eq:Tlimit4ModelC}
T \le T_{max} \left(\frac{1+e\cos{\theta}}{1+e}\right)^3
\end{equation}
That is, by imposing the condition given by \eqref{eq:ModelC} (i.e., Model $\widetilde{B}$) the maximum value of the spacecraft thrust is restricted by the right-hand-side of the inequality given in \eqref{eq:Tlimit4ModelC} (instead of $T_{max}$).  This artificial restriction of the spacecraft's actual thrust $T$ is simply an outcome of choice, namely, a selection of Model $\widetilde{B}$.  The reduced spacecraft thrust capacity is shown by the shaded region in Fig.~\ref{fig:ModelCspace} for $e = 0.20$.
%
%
It is apparent from Fig.~\ref{fig:ModelCspace} and \eqref{eq:Tlimit4ModelC} that the maximum value of the thrust capacity $T_{max}$ is indeed fully utilized by \eqref{eq:Tlimit4ModelC} for $\theta = 0$ (and $\theta = 2\pi$).  However, for most other values of $\theta$ (mod $2\pi$) Model $\widetilde{B}$ artificially limits the spacecraft thrust capacity to values well below $T_{max}$. In fact, it is clear from \eqref{eq:Tlimit4ModelC} that $T$ is most restricted at $\theta = \pi$.  Substituting $\theta = \pi$ in \eqref{eq:Tlimit4ModelC} we get,
\begin{equation}\label{eq:TlimitWorst}
T\big|_{\theta = \pi} \le T_{max} \left(\frac{1-e}{1+e}\right)^3
\end{equation}
A plot of the worst case limitation given by \eqref{eq:TlimitWorst} is shown in Fig.~\ref{fig:ModelCworstLimits}. For the Earth-Moon system (using $e = 0.05$), \eqref{eq:TlimitWorst} generates,
$$T\big|_{\theta = \pi} \le 0.74\,T_{max} $$
In other words, if Model $\widetilde{B}$ is used for spacecraft dynamics and control in cislunar space then this choice utilizes less than $75\%$ of the engine's thrust capability.

As a final point of comparison, consider the situation where a control solution $\theta \mapsto \widetilde{\bT}(\theta)$ is obtained by considering Model~$\widetilde{B}$. To implement this solution, the actual thrust $\bT$ on the spacecraft must be ``remodulated'' according to $\theta \mapsto (1+e\cos\theta)^3\,\widetilde{\bT}(\theta)$.  That is, even if the thrust program $\theta \mapsto \widetilde{\bT}(\theta)$ were to be a constant, the actual thrust on the spacecraft would require modulation.

\subsection{Introducing Model $\bC$: An Engine-Agnostic Dynamical Model}

All three models discussed so far employ the rocket mass flow rate equation to augment \eqref{eq: ERTBP_EOM_pseudoPot_ver}. A central issue in using \eqref{eq: ERTBP_timeMassFlowRate} is that a realistic value of $I_{sp}$ (typically paired with $T_{max}$) must be chosen for analysis\cite{ozimek_low-thrust_2010, oshima_global_2017, du_low-thrust_2023}. This limits the theoretical value of a particular computational analysis.  In fact, by choosing a very high value of $I_{sp}$, one can easily demonstrate a low value of propellant consumption. Hence, it is highly desirable to generate dynamical equations and a model for propellent consumption that is agnostic to the specifics of a particular engine.  The basic mathematics for this concept for a generic spacecraft was developed in  \cite{ross_how_2004,ross_space_2006}.  On the basis of these results and Model~$A$, we propose Model~$C$ for the dynamical equations of a thrusting spacecraft in an ER3BP according to:
\begin{subequations}\label{eq: ERTBP_EOM_pseudoPot_withUCntrl}
    \begin{align}
         x'' - 2y' &= \partial_x\Psi(x,y,z, \theta) + \frac{ u_x}{(1+e\cos{\theta})^3}\,
         \\
         y'' + 2x' &= \partial_y\Psi(x,y,z, \theta) + \frac{ u_y}{(1+e\cos{\theta})^3} \,
         \\
         z'' &= \partial_z\Psi(x,y,z, \theta) + \frac{ u_z}{(1+e\cos{\theta})^3} \,
    \end{align}
\end{subequations}
where, $u_x$, $u_y$ and $u_z$ are control variables.  Equation~(\ref{eq: ERTBP_EOM_pseudoPot_withUCntrl}) can also be written in a deceptively simple vector form as,
\begin{equation}
    \br'' + 2\,\hat{\bk} \times \br' = \nabla_{\br} \Psi(\br, \theta) + \frac{\bu}{(1+e\cos\theta)^3}
  \label{eq: ERTBP_pseudoPot_vecForm}
\end{equation}
where $\bu = (u_x, u_y, u_z)$ is a control vector that serve as a proxy for the thrust vector and is hence bounded by (see \eqref{eq:Tq-bound}),
\begin{equation}\label{eq:uq-bound}
\norm{\bu}_q \le u_{max}, \qquad \frac{1}{q} + \frac{1}{p} = 1
\end{equation}
As discussed in \cite{ross_how_2004,ross_space_2006} the correct model for propellant consumption is the $L^1$-norm of $t \mapsto \bu(t)$.  This \textit{does not} mean one must necessarily choose the (finite-dimensional) $\ell^1$-norm of $\bu \in \real{3}$.  In fact, from a purely mathematical perspective, one may choose any $\ell^p$-norm for $\bu$.  As shown in \cite{ross_how_2004,ross_space_2006}, just three values of $p \in \set{1, 2, \infty}$ capture the basic ideas of thruster configurations with $p = 1, 2$ being the most common.  In any case, the model for propellant consumption proposed in \cite{ross_how_2004,ross_space_2006} is given by the $\ell^p$ variant of the $L^1$ norm of $t \mapsto \bu(t)$,
\begin{equation}\label{eq:proxByDef}
\Delta prox:=  \int_{t_0}^{t_f} \|\bu(t)\|_p  \,dt
\end{equation}
In \eqref{eq:proxByDef} we use the symbol $\Delta prox$ to denote the fact that the the $L^1$-norm does not compute the actual propellant, rather a proxy that accurately measures the $I_{sp}$-agnostic propellant consumption. To clarify this last statement, we first compute the actual amount of propellant consumption, $\Delta prop$, by integrating \eqref{eq: ERTBP_trueAnomMassFlowRate_withCharTime}
\begin{multline}
    \Delta prop = m_{sc}(\theta_0) - m_{sc}(\theta_f) = - \int_{\theta_0}^{\theta_f} m_{sc}' \, d\theta \\
    = \frac{1}{I_{sp}\,\omega_C \,g_0} \int_{\theta_0}^{\theta_f} \frac{\|\bT(\theta)\|_p}{(1+e\cos{\theta})^2} \,d\theta
    \label{eq: ERTBP_massChangeIntegral}
\end{multline}
Changing the domain of integration in \eqref{eq:proxByDef} from $t$ to $\theta$ we get,
\begin{equation}\label{eq:changeDom}
\Delta prox:=  \int_{\theta_0}^{\theta_f} \|\bu(\theta)\|_p  \left(\frac{dt}{d\theta}\right)\, d\theta  = \int_{\theta_0}^{\theta_f} \frac{\|\bu(\theta)\|_p}{\omega} \,d\theta
\end{equation}
Substituting \eqref{eq:omega=funtheta} in \eqref{eq:changeDom} and scaling the result by $\omega_C$, we get
\begin{equation}
    \Delta prox=  \int_{\theta_0}^{\theta_f} \frac{\|\bu(\theta)\|_p}{(1+e\cos{\theta})^2} \,d\theta
    \label{eq:DeltaProx4ER3BP}
\end{equation}
Comparing \eqref{eq:DeltaProx4ER3BP} and \eqref{eq: ERTBP_massChangeIntegral} it follows that $\Delta prox$ is an accurate engine-agnostic measure of $\Delta prop$. The propellant consumption measured by $\Delta prox$ has the (nondimensional) units of ``Delta-$V$'' that automatically accounts for ``finite-burn losses''\cite{biblarz,robbins-deltaV}.
If  $p=1$, then \eqref{eq:DeltaProx4ER3BP} is the $\ell^1$-$L^1$-norm corresponding to propellant consumption associated with six engines placed orthogonally along the $x$, $y$, $z$ axes\cite{ross_how_2004}.  If $p=2$, then \eqref{eq:DeltaProx4ER3BP} is the $\ell^2$-$L^1$-norm corresponding to propellant consumption associated with a single engine\cite{ross_how_2004}.  See \cite{ross_how_2004,ross_space_2006} for further details.  In summary, Eqs.~(\ref{eq: ERTBP_EOM_pseudoPot_withUCntrl}), (\ref{eq:uq-bound}) and (\ref{eq:DeltaProx4ER3BP}) constitute Model~$C$ for the control constrained dynamical equations for the ER3BP together with an engine-agnostic model for propellant consumption.

Per the comments in Section~\ref{sec:intro}, we briefly note that although the $\ell^p$-norm is not differentiable, it can easily be mapped to continuously differentiable functions through the addition of new control variables\cite{ross-book}. In other words, it is not necessary to use quadratic cost functionals on the basis of differentiability. More importantly, quadratic cost functionals generate non-fuel-optimal solutions with propellant consumption that can be as high as $50\%$ in excess of the optimal solution\cite{ross_space_2006}.  This problem is easily avoided using elementary control variable transformations\cite{ross-book}.

Lastly, we briefly note that we can define two other engine-agnostic models along the lines of \eqref{eq:Ttilde2T}.
Per the discussions in the previous subsections, these models will inherit all the issues noted therein; hence, we will not discuss them further.

\section{A Modified Formula for Impulsive $\Delta V$ Computation}
Mathematically, the quantity $\Delta prox$ defined in \eqref{eq:proxByDef} is simply the $L^1$ norm of the $\ell^p$ norm of $\bu(t)$. The physical units of $\Delta prox$ is the same as velocity but it is not the same as ``Delta-$V$'' as widely understood\cite{biblarz,hale}. For instance, $\Delta prox$ incorporates ``gravity loss'' and other loss terms (like drag loss, under atmospheric effects) whereas $\Delta V$ does not\cite{robbins-deltaV}.

\subsection{An ER3BP-Modification to the Rocket Equation}
There are multiple ways to derive the rocket equation.  Here, we take the textbook approach that is equivalent to dropping all terms in \eqref{eq: ERTBP_EOM_pseudoPot_ver_withNonDmass} except for the thrust and second derivative terms. Performing this exercise and writing the result in scalar form we get,
\begin{subequations}
\begin{align}
v' &= \frac{T/(m_C\, A_C)}{m(1+e\cos\theta)^3} \label{eq:vprime4rocket}\\
m' & = -\frac{T/(m_C\, g_0)}{(I_{sp}\omega_C)(1+e\cos\theta)^2} \label{eq:mprime4rocket}
\end{align}
\end{subequations}
where $v' = r''$.
Dividing \eqref{eq:vprime4rocket} by \eqref{eq:mprime4rocket} we get,
\begin{equation}\label{eq:dvbydm}
\frac{dv}{dm} = -\frac{(I_{sp} g_0) (\omega_C/A_C)}{m  (1+e\cos\theta)}
\end{equation}
Integrating \eqref{eq:dvbydm} under the standard $\Delta V$ assumption, namely, that the burn duration is infinitesimally short, we get,
\begin{equation}\label{eq:TaDa}
\Delta m = m_0\left[1 - \exp\left(-\frac{\Delta V (1 + e \cos\theta_I)}{(I_{sp}g_0)/v_C} \right) \right]
\end{equation}
where $\Delta m = m(\theta_0) - m_f (\theta_f)$ is the change in mass, $\Delta V = v(\theta_f) - v(\theta_0)$ is the corresponding change in velocity, $m_0 = m(\theta_0)$, $v_C = A_C/\omega_C$ (Cf.~\eqref{eq:vC=})  and $\theta_I$ ($= \theta_f = \theta_0$) is the point of application of the impulse.  Note that $\Delta V$ in \eqref{eq:TaDa} is not in physical units but consistent with the (nondimensional) velocity units implied by $r' = v$.

\subsection{Discussion of \eqref{eq:TaDa}}
Perhaps, the most glaring modification to the rocket equation that is apparent in \eqref{eq:TaDa} is its dependence on the true anomaly of the primaries. This dependence vanishes if $e=0$, i.e., a CR3BP. In this case, the nondimensional $\Delta V$ is simply multiplied by the circular speed, $v_C$ (of the primaries) and \eqref{eq:TaDa} reduces to Tsiolkovsky's famous formula.  This point once again demonstrates in a different way how the CR3BP hides the nuance implied in \eqref{eq:TaDa}.

When $e \ne 0$, the numerator in the exponent of \eqref{eq:TaDa} is attenuated or amplified by an amount that depends upon the point of application of the $\Delta V$ maneuver. The largest attenuation occurs at the apoapse ($\theta_I = \pi$) while the largest amplification is produced at the periapse ($\theta_I = 0$). Obviously, neither attenuation nor amplification is incurred if the impulse is applied at the semilatus rectum ($\theta_I = \pi/2$).

If a propellent-optimal maneuver is designed in an ER3BP based on impulsive maneuvers considerations, then the mission Delta-$V$ cost, $(\Delta V)_{mission}$ that correctly represents the propellant consumption according to \eqref{eq:TaDa} is given by,
\begin{equation}\label{eq:TaDa-total}
(\Delta V)_{mission} = \sum_{i=1}^N \Delta V_i (1 + e \cos\theta_i)
\end{equation}
where $\theta_i, i = 1, \ldots, N$ are the individual locations of the $N$ impulses, $\Delta V_i, i = 1, \ldots, N$.  In other words, unlike the CR3BP, the computation of the mission $\Delta V$ is not agnostic to the point of application of the individual impulsive maneuvers.  Thus, a mission design in the ER3BP that uses \eqref{eq:TaDa-total} with $e=0$ would have either (a) overestimated or underestimated the total Delta-V and/or (b) produced a solution with incorrect values of $\Delta V_i, i = 1, \ldots, N$ including the number, location and values of the impulses.

When considered in isolation, \eqref{eq:TaDa} (or \eqref{eq:TaDa-total}) appears to suggest that the propellant-optimal location of an impulse is the apoapse point $\theta_I = \pi$. Although the apoapse location provides the largest attenuation of the effective $\Delta V$, it may not be optimal location from the point of view of a particular mission when the entire suite of boundary conditions are taken into account.  In other words, it is entirely possible for a large number of impulses to be concentrated at the ``worst points,'' $\theta_i = 0, i = 1, \ldots $ in order to satisfy boundary conditions. Conversely one may easily argue the opposite of the preceding argument by suggesting that the $\Delta V$ at the apoapse point must necessarily be larger than the one at the periapse point to counter the attenuation produced at the former location.  Needless to say, \eqref{eq:TaDa-total} must be viewed holistically in terms of the mission design and not in isolation.

\subsection{Connections Between $\Delta prop$, $\Delta prox$ and $\Delta V$}
Obviously, the ``exact'' measure of propellant consumption is $\Delta prop$ given by \eqref{eq: ERTBP_massChangeIntegral}.  However, as noted earlier and elsewhere\cite{ross_how_2004,ross_space_2006} because \eqref{eq: ERTBP_massChangeIntegral} is $I_{sp}$-dependent, there is a need to design and analyze missions before a specific propulsion system is selected.  The rocket equation and $\Delta V$ minimization has been the core of such analyses since the dawn of the space age. Because a practical propulsion system generates finite burns, there has been a substantial body of research focused on minimizing ``losses'' incurred by gravity, drag etc. These ideas are predicated on the assumption that the optimal location of the finite burns would be ``centered'' at the point of application of the Delta-$V$'s. It was suggested in [\citenum{ross_how_2004}] and [\citenum{ross_space_2006}] that taking the $L^1$ functional would obviate the need to determine the loss terms while simultaneously ``relocating'' the finite burns without using any $I_{sp}$-specific information. The standard $L^1$ functional model holds for any space mission where the independent variable is time. In the ER3BP dynamical model (Model $C$), the independent variable is true anomaly and not time. If the true-anomaly to time transformation affected the dynamical equations and the mass-flow rate equations the same way, then all of the results developed in [\citenum{ross_how_2004}] and [\citenum{ross_space_2006}] would continue to hold without change. As it turns out, Nechvile's transformation affects the mass flow rate equation differently than the dynamical equations.  The latter contains a cubic term (Cf.~\eqref{eq:cubic}) while the former is impacted by a quadratic (Cf.~\eqref{eq:quadratic}).  The ratio of these two terms is ``linear'' which appears in \eqref{eq:TaDa} (and hence \eqref{eq:TaDa-total}).

In \eqref{eq:DeltaProx4ER3BP}, the $L^1$ measure of propellant consumption is modified by the quadratic cosine term.  This formula appears to be at loggerheads with \eqref{eq:TaDa} which presents a modification to the rocket equation by a linear cosine term. To prove that both measures are consistent we use the notion embedded in \eqref{eq:TaDa-total}.  Replacing the summation in \eqref{eq:TaDa-total} by an integral, we get,
\begin{equation}\label{eq:TaDa-integral}
(\Delta V)_{mission} = \int_{v_0}^{v_f} (1 + e \cos\theta)\, dv = \int_{\theta_0}^{\theta_f} (1 + e \cos\theta) v' \, d\theta
\end{equation}
Substituting \eqref{eq:vprime4rocket} on the right-hand-side of the second equality in \eqref{eq:TaDa-integral} recovers the quadratic cosine term of \eqref{eq: ERTBP_massChangeIntegral} which was the basis for defining \eqref{eq:DeltaProx4ER3BP}.

\section{An Illustrative Cislunar Space Mission Scenario}\label{sec:missionEx}
The many number of nuances and generic mathematical results derived and discussed in the preceding sections can be numerically illustrated by considering a sample mission scenario in the ER3BP.  To this end, we consider the problem of transferring a spacecraft from a 2:1 $L_1$-Lyapunov orbit to a 4:1 $L_2$ near rectilinear Halo orbit (NRHO) in the Earth-Moon ER3BP with $e=0.0549$. Complete mathematical details of this problem formulation are described in \cite{ER3BP-Hawaii-2025}.  Here, we limit the discussions to just the use and impact of of different cost functionals for propellant consumption.

To map the mathematics of the prior sections to the physics of the posed problem, we assume the spacecraft is equipped with six equal thrusters, one on each side of its face and labeled $u_i^+, u_i^-, i = 1, 2, 3$.  This is the $\ell^1$ configuration described in [\citenum{ross_how_2004}] and [\citenum{ross_space_2006}].   We consider an optimization of trajectories for the same mission (i.e. boundary conditions) generated by three different cost measures given by,
\begin{itemize}
\item[a)] the quadratic functional,
\begin{equation}\label{eq:JQ=}
    J_Q :=  \int_{\theta_0}^{\theta_f}\sum_{i=1}^3 \left(\left[u_i^+(\theta)\right]^2 + \left[u_i^-(\theta)\right]^2\right) \,d\theta
\end{equation}
\item[b)] an $L^1$-functional without the cosine term,
\begin{equation}\label{eq:J1=}
    J_1 :=  \int_{\theta_0}^{\theta_f}\sum_{i=1}^3 \left(u_i^+(\theta) + u_i^-(\theta)\right) \,d\theta
\end{equation}
and
\item[c)] an  $L^1$-functional with the inclusion of the cosine term,
\begin{equation}\label{eq:Dprox=}
     \Delta prox :=  \int_{\theta_0}^{\theta_f} \frac{\sum_{i=1}^3 \left( u_i^+(\theta) + u_i^-(\theta)\right)}{(1 + e\cos\theta)^2} \,d\theta
\end{equation}
\end{itemize}
In all cases, all of the six thrusters are bounded by $[0, 0.1]$.  The numerical optimization is performed using an $\alpha$-version of the software package DIDO\textsuperscript{\textcopyright}\cite{DIDO:arXiv}.  This variant of DIDO implements a guess-free version\cite{auto-knots,ross:guess-free,spec-alg} of the universal Birkhoff theory\cite{newBirk-part-I,newBirk-part-II,newBirk-bvp,ross:Hessians} for trajectory optimization.

\begin{figure}[h!]
\centering
    \includegraphics[width=\columnwidth, trim={1.75in, 0, 2in, 0}, clip]{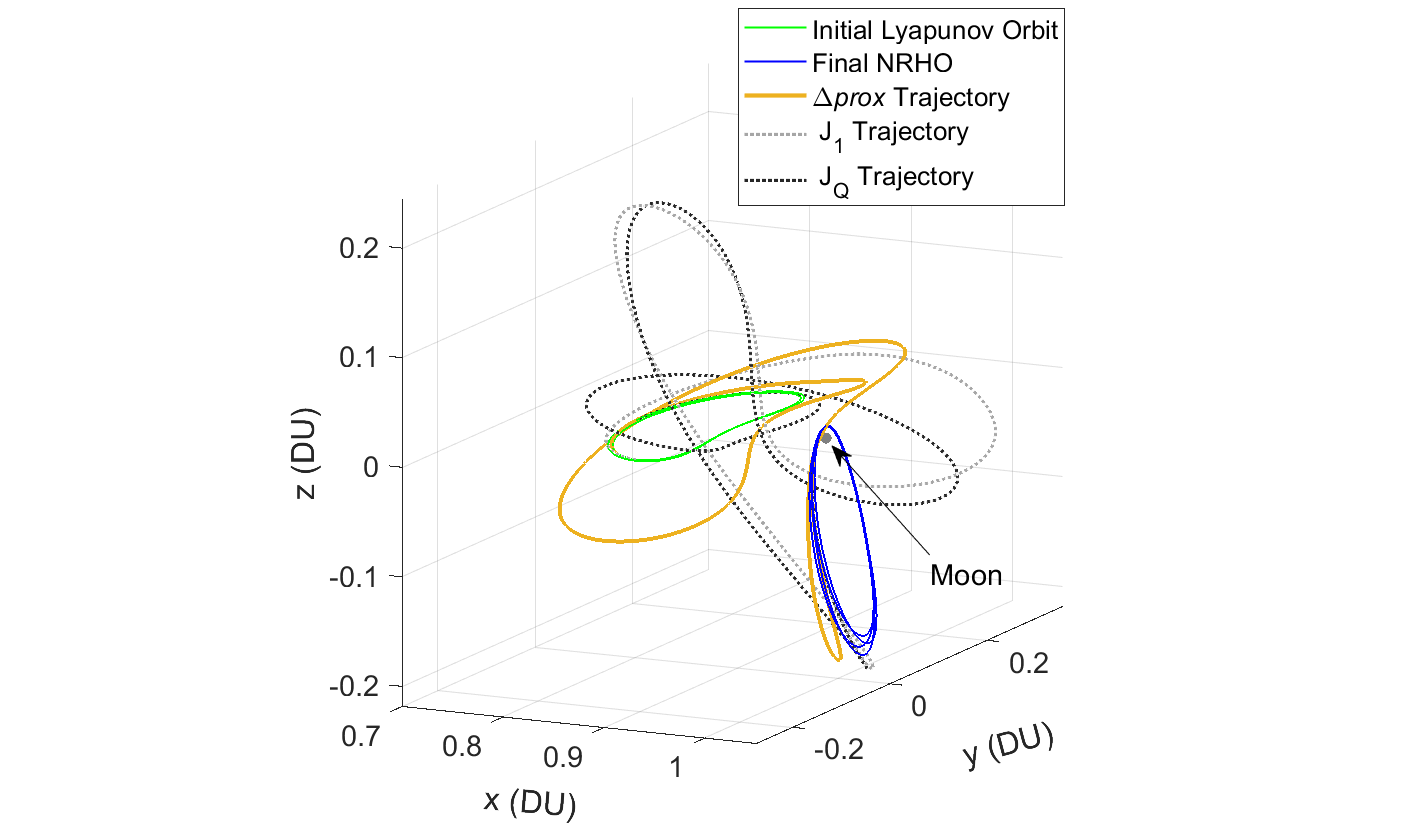}
\caption{Three trajectories in the cislunar Nechvile frame generated by minimizing Eqs.~(\ref{eq:JQ=}), (\ref{eq:J1=}) and (\ref{eq:Dprox=}).}
\label{fig:three-trajs}
\end{figure}

The computed optimal trajectories in the Nechvile frame are displayed in Fig.~\ref{fig:three-trajs} for all three cases.
%
%
These numerically optimized trajectories illustrate the following:
\begin{enumerate}
\item All three trajectories are different.
\item Although the $J_1$ and $\Delta prox$ cost functionals differ only by the quadratic cosine term (Cf.~\eqref{eq:J1=} and \eqref{eq:Dprox=}), the resulting two trajectories are substantially different from one another.
\item The $J_1$- and $J_Q$-trajectories appear close together although their cost functionals are substantially different.  The only thing common between the $J_Q$ and $J_1$ functionals is that they do not incorporate the quadratic cosine term. Both trajectories use incorrect measures of propellant consumption.
\end{enumerate}
The amount of propellant consumed as given by the $\Delta prox$ measure for all three trajectories are provided in  Table~\ref{tab:prop}.  The column marked Propellant Consumed shows that using $\Delta prox$ as a cost functional generates the lowest amount of propellant consumption as it should\cite{ross_how_2004,ross_space_2006}.  The $J_Q$ cost functional generates the highest amount of propellant consumption and consistent with the theoretical discussions of Refs~[\citenum{ross_how_2004}] and \cite{ross_space_2006}. As shown by the column marked \% Excess in Table~\ref{tab:prop}, the penalty for ignoring the quadratic cosine term, namely, the use of the $J_1$ cost functional, is 81\%.  The penalty for using the quadratic cost is higher by more than a factor of two (i.e., 114\%).
%
\begin{table}[h!]
    \centering
    \caption{Propellant Consumption Incurred by Different Cost Functionals for the Cislunar Mission Shown in Fig.~\ref{fig:three-trajs}.}
    \label{tab:prop}
    \begin{tabular}{lccc}
        \toprule
        \textbf{Cost Functional} & \multicolumn{2}{c}{\textbf{Propellant Consumed}} & \multicolumn{1}{c}{\textbf{Value of the}} \\
        \cmidrule(lr){2-3} 
        & $\mathbf{\Delta prox}$ & \textbf{\% Excess} &  \textbf{Cost Functional}  \\
        \midrule
        $J_Q$ &         0.702 & 114 &  0.176\textsuperscript{a} \\
        $J_1$ &         0.594 & 81 &   0.545 \\
        $\Delta prox$ & 0.328 & 0 &   0.328 \\
        \bottomrule\\[-0.5em]
        \multicolumn{4}{l}{$^a$ \footnotesize{To make it consistent with the ``units'' of $J_1$ and $\Delta prox$, the value listed }}\\
        \multicolumn{4}{l}{\footnotesize{\ here is $\sqrt{J_Q}$.}}
    \end{tabular}
    \vspace{2mm}
\end{table}
%

It is important to note that the numbers listed in Table~\ref{tab:prop} are valid only for the trajectories shown in Fig.~\ref{fig:three-trajs}.  That is, these numbers will most likely be different (lower or higher) for different mission scenarios.  What is theoretically provable\cite{ross_how_2004,ross_space_2006} is that the propellant consumption for the $J_1$- and $J_Q$-optimized trajectories will always be greater than (or, at best, equal to) the one obtained by minimizing $\Delta prox$.

The numbers in the last column of Table~\ref{tab:prop} are quite interesting and may possibly mislead an analyst. First, we note that the optimized value of $J_Q$ was found to be $0.031$.  Second, because $J_Q$ is quadratic, the value of $\sqrt{J_Q}$ is provided in Table~\ref{tab:prop} so that it can be ``compared'' to the remainder of the numbers in the table.  From the perspective of the last column of Table~\ref{tab:prop}, an analyst might erroneously conclude that the $J_Q$ generated trajectory is the ``best'' even though it is the worst performer in terms of propellant consumption. Another point of misdirection is the difference between the $J_1$ cost and the $J_1$ propellant consumption listed in the second row of Table~\ref{tab:prop}.  The difference between these two numbers is approximately 9\%.  This difference is due to the exclusion or inclusion of the quadratic cosine term (compare \eqref{eq:J1=} and \eqref{eq:Dprox=}).  In the absence of the third row, one might erroneously conclude that the quadratic cosine term makes only a ``small'' difference in the computation of the propellant consumption (wherein small implies less than 10\%). That this conclusion is misleading is apparent from the values of the third row of Table~\ref{tab:prop}.  In other words, when the cost functional is different by a presumable ``small correction'' term given by the inverse of $(1 + e\cos\theta)^2$, the resulting trajectory can be substantially different in generating a much lower propellant solution (than a mere difference of 9\%).  In this example mission scenario, the lower propellant solution saves 81\% of the fuel consumed.

\begin{figure}[h!]
\centering
    \includegraphics[width=\columnwidth]{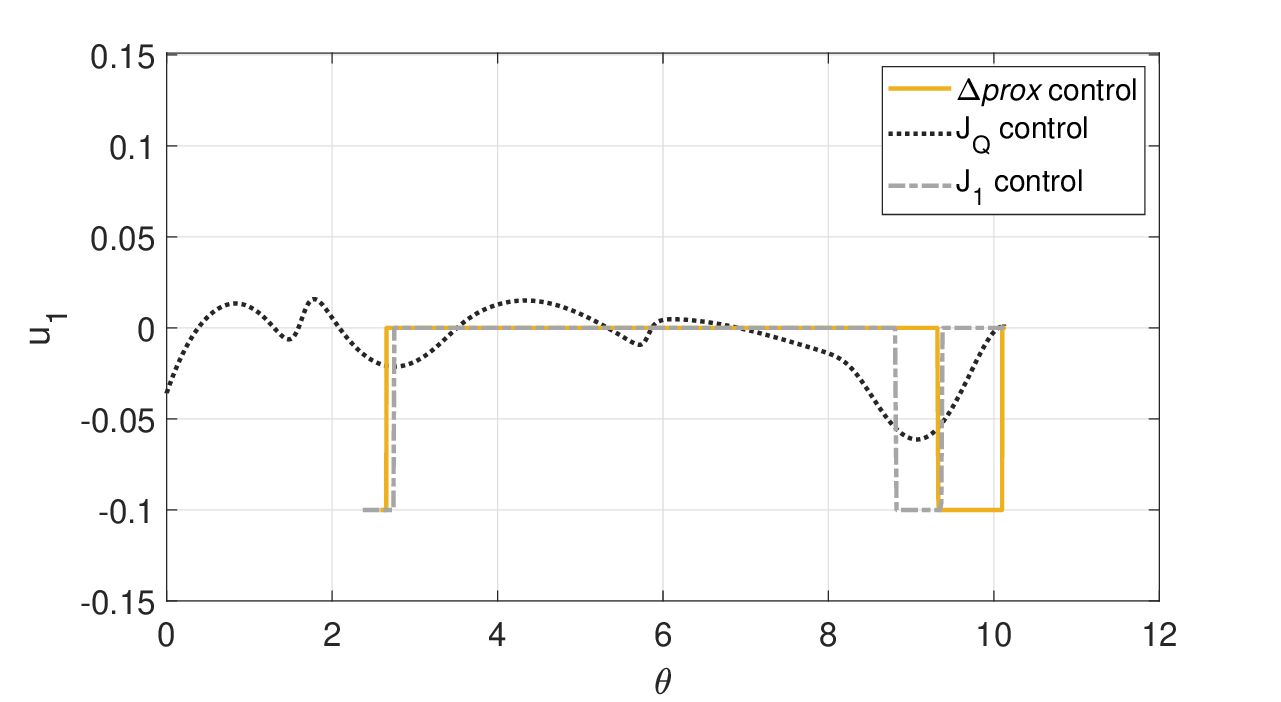}
    \includegraphics[width=\columnwidth]{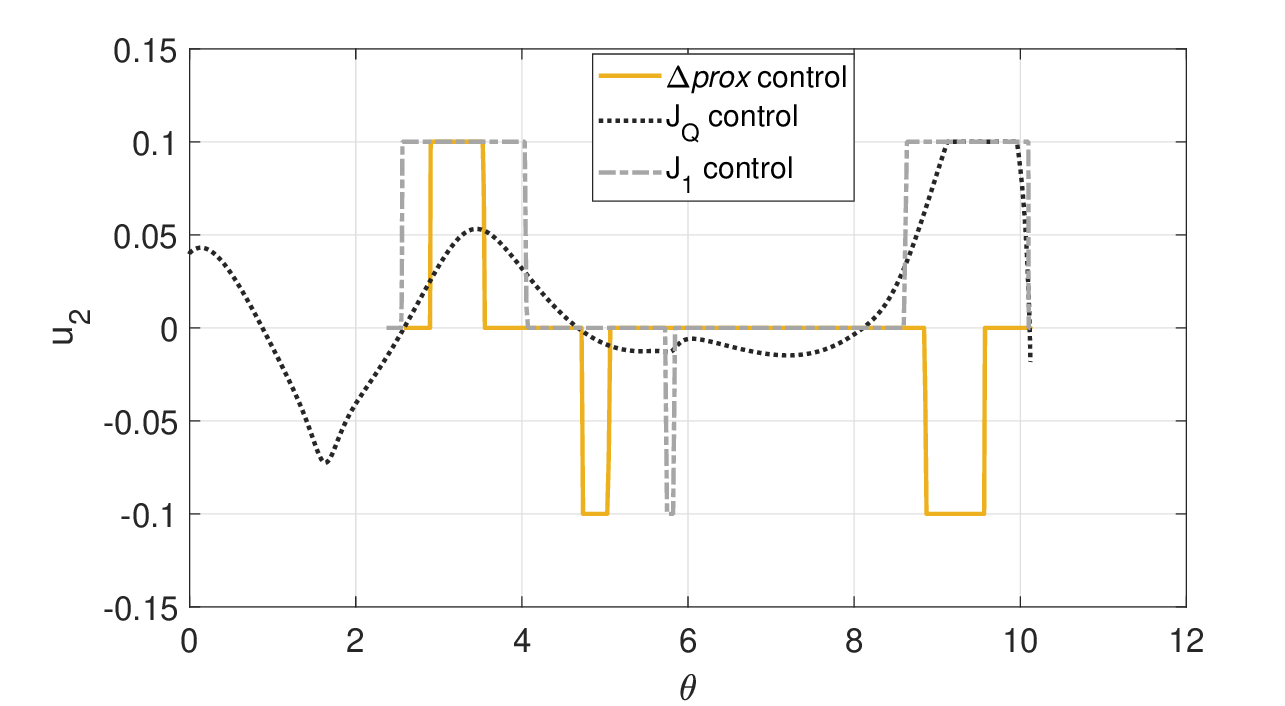}
    \includegraphics[width=\columnwidth]{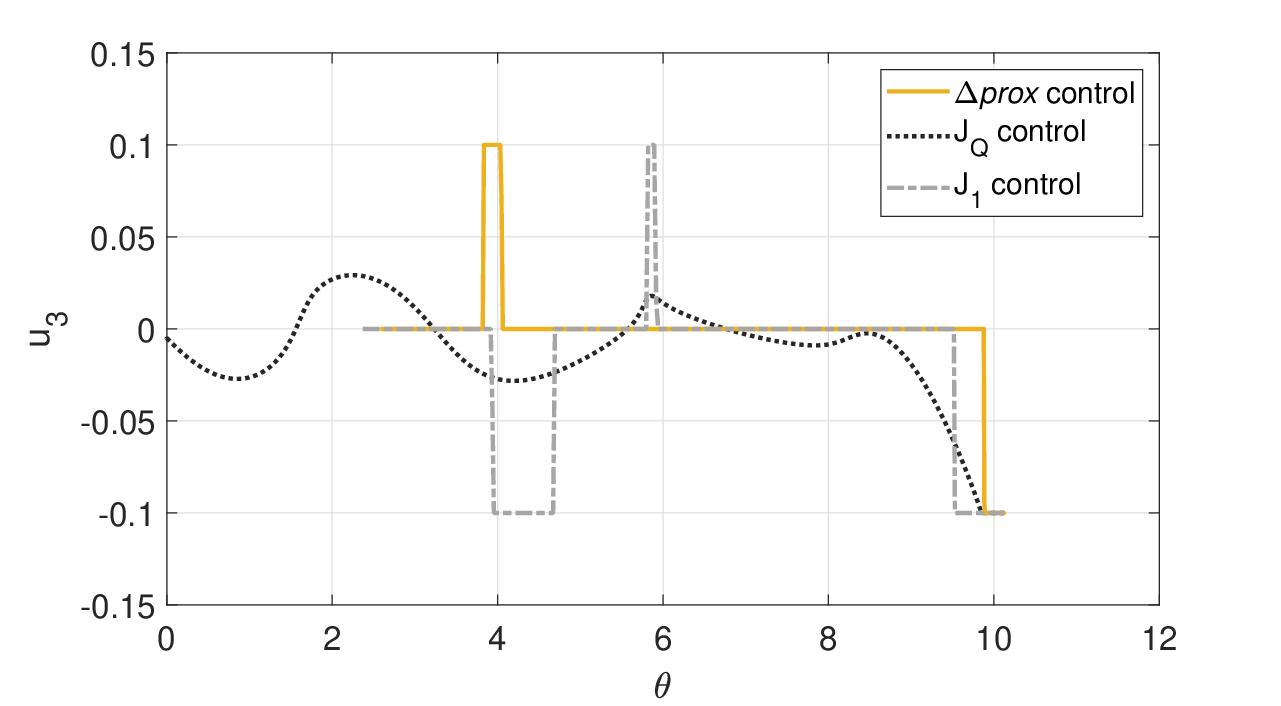}
\caption{Control profiles for the three trajectories shown in Fig.~\ref{fig:three-trajs}.}
\label{fig:three-controls}
\end{figure}

The differences between the trajectories shown in Fig.~\ref{fig:three-trajs} and the stark contrast in numbers displayed in Table~\ref{tab:prop} are more vivid when one examines the control solutions for all three cases. The components of the the control trajectories $\theta \mapsto \bu(\theta)$ for all three cases are shown in Fig.~\ref{fig:three-controls}.
%
%
These components correspond to $\bu(\theta) := (u_1(\theta), u_2(\theta), u_3(\theta))$, where $u_i(\theta) = u_i^+ - u_i^-, i = 1, 2, 3$.  As a result, $u_i \in [-0.1, 0.1],\ i = 1, 2, 3$.  These control plots showcase the following:
\begin{enumerate}
\item As expected\cite{ross_how_2004,ross_space_2006}, $\theta \mapsto \bu(\theta)$ is continuously modulated for the trajectory generated by the $J_Q$ cost functional.
\item The $J_1$-control solution is substantially different from that of $\Delta prox$ even though they are both $L^1$ functionals, albeit the former being an incorrect measure of propellant consumption.
\item The $J_1$ control solution is quite different from that of the $J_Q$ solution even though they both generate (topologically) similar three-dimensional trajectories as shown in Fig.~\ref{fig:three-trajs}.
\item The $J_1$ solution may also be interpreted in terms of generating finite-burn equivalents of impulsive solutions using \eqref{eq:TaDa-total} but without the linear cosine correction term. Similarly, the $\Delta prox$ solution is the finite-burn equivalent of \eqref{eq:TaDa-total}. Under these interpretations, it is obvious from the plots shown in Fig.~\ref{fig:three-controls} that the number, location and magnitude of the ``Delta-$V$'s'' of the $J_1$ solution are substantially different from the one generated by minimizing $\Delta prox$. For example, the component of ``impulse'' in the \#3 direction (see $u_3$ plot in Fig.~\ref{fig:three-controls}) contains an extra impulse for the $J_1$ functional minimization (when compared with the $\Delta prox$ solution).
\end{enumerate}
In summary, this example shows that small eccentricities in the restricted three body problem can have big impacts in trajectory optimization.

\section{Conclusions}\label{sec:conclusions}
The addition of the mass-flow rate equation to Nechvile's equations for the ER3BP generates an interesting phenomenology that can be directly linked to pulsation. The pulsation vanishes for the CR3BP; hence, intuitive modifications, or the lack of it, when considering circular cases do not transfer to elliptic situations.  Even when the eccentricity is small, as in the case of the Earth-Moon system, the ensuing results from trajectory optimization may be dramatically different. In principle, the dramatic differences in optimized trajectories is predictable from the fact that the restricted three-body problem exhibits ``pockets of chaos.'' Consequently, it should not be a surprise that results from a CR3BP do not carry over to the elliptic case with ``small modifications for small eccentricities'' even with the inclusion of corrective guidance maneuvers. When these ideas are adjoined to the use of incorrect measures of propellant consumption, the ensuing results are amplified even more.  Mathematically, the main culprit for the apparent disruption of well-trodden ideas is the cubic-quadratic mismatch in the thrust terms in the ER3BP.  Taking this mismatch into account leads to a modification of the famous  Tsiolkovsky formula.  For finite thrust arcs, the $L^1$-functional that correctly captures the propellant consumption also undergoes a modification due to Nechvile's transformation.  Results from a sample mission scenario shows that ignoring these apparent nuances can lead to starkly different trajectories and propellant consumptions. In other words, it may be necessary even in the early stages of mission design considerations to not ignore small eccentricities because their impact may be quite large.  All of these conclusions hold for cislunar space trajectories because the eccentricity of the Earth-Moon system is not zero. The non-zero value of the Earth-Moon eccentricity may be harnessed for lowering propellant consumption via new pathways.

\section*{Acknowledgments}
Funding for this research was provided by the Air Force Office of Scientific Research (AFOSR) and the U.S. Navy.  The views and conclusions contained herein are those of the authors and should not be interpreted as necessarily representing the official policies or endorsements, either expressed or implied, of the AFOSR or the U.S. Navy or the U.S. Government.

\end{document}